\documentclass{article}
\usepackage{amssymb, amsmath}
\usepackage[latin9]{inputenc}
\usepackage[francais]{babel}
\usepackage[affil-it]{authblk}
\usepackage{graphicx}
\usepackage{hyperref}
\usepackage{color}
\usepackage[onehalfspacing]{setspace}
\usepackage[nottoc, notlof, notlot]{tocbibind}
\usepackage{ulem}
\makeatletter
\def\imod#1{\allowbreak\mkern10mu({\operator@font mod}\,\,#1)}
\makeatother

\newcommand{\Q}{\mathbb{Q}}

\begin{document}

\title{Un portrait kaléidoscopique \\ du jeune Camille Jordan}
\author{Fr{\'e}d{\'e}ric Brechenmacher%
\thanks{Electronic address: \texttt{frederic.brechenmacher@euler.univ-artois.fr} \\ \textit{Ce travail a bénéficié d'une aide de l'Agence Nationale de la Recherche : projet CaaF{\'E} (ANR-10-JCJC 0101)}}}
\affil{Université d'Artois, \\ Laboratoire de math{\'e}matiques de Lens (EA 2462) \\ rue Jean Souvraz S.P. 18, F- 62300 Lens France}
 \date{}
\maketitle

\begin{center}
\includegraphics[scale=0.9]{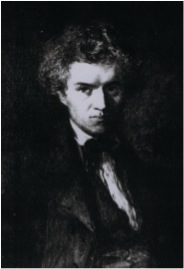}\\*
Camille Jordan (1838-1922), sous le pinceau de son oncle, \\ le peintre Puvis de Chavanne.
\end{center}

\tableofcontents

\newpage

Un premier visage de Camille Jordan se dessine ci-dessus à l'âge de 17 ans auquel notre héros intègre l'École polytechnique (1855). Cet article est principalement consacré à la première décennie de travaux de ce mathématicien, c'est-à-dire à la période séparant la soutenance en 1860 par Jordan de sa première thèse d'algèbre et la parution en 1870 du célèbre \textit{Traité des substitutions et des équations algébriques}, second héros de cet article.\\

\begin{center}
\includegraphics[scale=0.6]{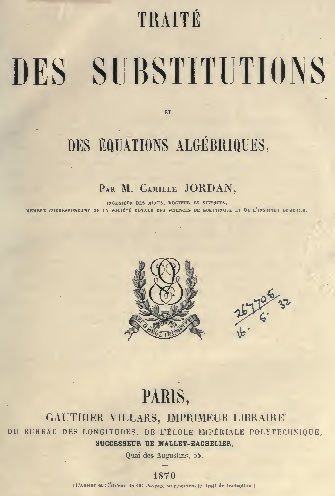}\\*
\end{center}

Comme nous allons le voir, le visage du jeune Jordan nous présente différentes facettes qu'il nous faudra scruter avec attention afin de restituer l'originalité des approches développées par le mathématicien dans son \textit{Traité}.

D'un point de vue institutionnel, Jordan poursuit une carrière classique de savant du XIX\up{e} siècle : issu de la bourgeoisie lyonnaise, il intègre l'École polytechnique puis débute ses travaux mathématiques en parallèle de ses fonctions d'ingénieurs des mines.\footnote{Pour une notice biographique sur Jordan, voir notamment : LEBESGUE (Henri), Notices d'histoire des mathématiques. Notice sur la vie et les travaux de Camille Jordan, \textit{L'enseignement mathématique} (1923), p. 40-49.}

A partir de la fin des années 1870, Jordan se voit confié les principales clés des institutions mathématiques en France : nommé Professeur d'Analyse à l'École polytechnique en 1876, élu à l'Académie des sciences en 1881, puis nommé Professeur au Collège de France en 1883. 

Jordan est également célèbre pour son \textit{Cours d'Analyse de l'École polytechnique} dont le premier volume paraît en 1882.\footnote{Voir à ce sujet GISPERT (Hélène), \textit{Camille Jordan et les fondements de l'analyse : Comparaison de la 1ère édition (1882-1887) et de la 2ème (1893) de son cours d'analyse de l'Ecole polytechnique},  Thèse de doctorat, Orsay : 1982.} 

\begin{center}
\includegraphics[scale=0.5]{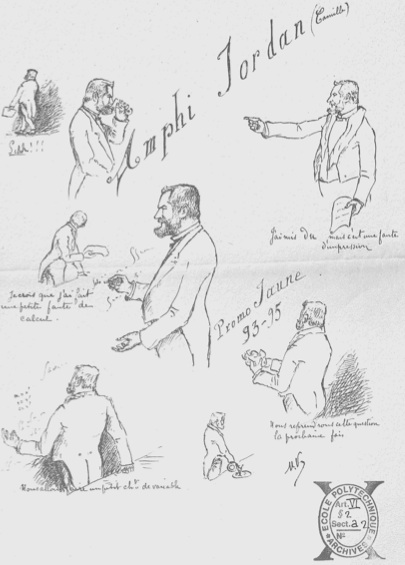}\\*
Archives de l'École polytechnique.
\end{center}

On peut l'apercevoir ci-dessus, croqué dans la force de l'âge par un élève de l'École polytechnique : traçant à la craie des successions d'intégrales ; caricaturé reproduisant une faute d'impression de son Cours d'Analyse : "je crois que j'ai fait une petite faute de calcul" ; prenant le temps de la réflexion en vidant un verre d'eau sucrée : "nous allons faire un petit changement de variable" ; déambulant sur l'estrade ; puis, chiffon à la main : "nous reprendrons cette question la prochaine fois". Jordan était visiblement apprécié des élèves polytechniciens qui ont baptisé de son nom le verre d'eau sucrée du professeur.

\begin{center}
\includegraphics[scale=0.5]{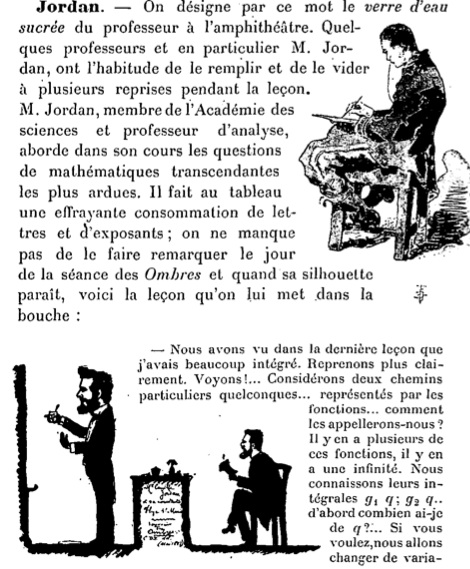}
\includegraphics[scale=0.5]{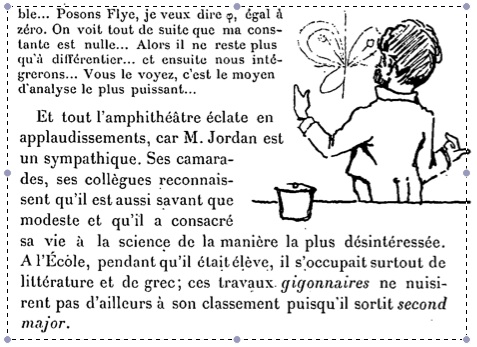}\\*
Testard, É. (ed). \textit{L'Argot de L'X}, \\ Paris: Albert-Lévy et G. Pinet, 1894, p. 178.
\end{center}

De 1885 à sa mort en 1922, Jordan était également directeur du \textit{Journal de mathématiques pures et appliquées}, l'un des principaux journaux de recherche mathématique de l'époque.\footnote{Voir BRECHENMACHER (Frédéric), \href{http://sabix.revues.org/730}{Le  "journal de M. Liouville" sous la direction de Camille Jordan (1885-1922)}, \textit{Bulletin de la Sabix}, 45 (2009), p. 65-71.} 

A première vue, notre héros semble donc parfaitement en prise avec l'organisation institutionnelle des sciences mathématiques de son temps. Pourtant, Jordan a souvent été considéré comme "travaillant dans une solitude presque totale",\cite{9} isolé sur la scène mathématique française et même "presque allemand".\footnote{KLEIN (Felix), [1921-1923] \textit{Gesammelte mathematische Abhandlungen}, Springer: Berlin, 1921, vol. 1, p. 51} 

Brosser un portrait de Camille Jordan pose donc le problème de la restitution des dimensions collectives dans lesquelles saisir la création mathématique individuelle. Jordan était-il réellement isolé ? Dans le cas contraire, dans quels cadres collectifs ses travaux se situaient-ils ? 

Afin d'aborder ce problème, il est nécessaire d'étudier à la fois les espaces sociaux dans lesquels se situent les mathématiques et les organisations collectives des savoirs en théories ou disciplines. Comme nous allons le voir, la présentation de Jordan comme savant isolé participe d'un portrait qui a souvent forcé le trait sur le développement de la "théorie des groupes", de la "théorie de Galois" et plus généralement de l'"algèbre".  Or aucune de ces catégories n'est stable historiquement ; il faut donc en restituer les significations variées dans différents temps et espaces. 

Nous allons envisager ces questions en prenant à la fois pour héros un mathématicien, Jordan, et un texte mathématique, le \textit{Traité des \textcolor{blue}{substitutions}\footnote{\textcolor{blue}{Une substitution est ce que nous désignerions aujourd'hui comme une permutation d'un nombre fini de lettres ($a$, $b$, $c$, $d$, ...). En des termes qui nous sont contemporains, la théorie des substitutions correspond donc à la théorie des groupes finis.}} et des équations algébriques}. \\

Ce texte a été conçu afin de permettre différents niveaux de lectures : des compléments mathématiques peuvent être obtenus en cliquant sur les termes colorés en \textcolor{blue}{bleu}, des citations peuvent être déroulées en cliquant sur les termes colorés en \textcolor{red}{rouge}, enfin des compléments historiques sont proposés dans quatre annexes accessibles par des liens hypertextes.

\section{Un "grand algébriste" en relation directe avec les "idées de Galois"} 

Un premier portrait du jeune Camille Jordan s'esquisse en clair obscur, dans l'ombre de la grande figure d'un autre jeune mathématicien: Évariste Galois. Les nombreuses nécrologies publiées après le décès de Jordan en 1922 insistent en effet toutes sur la relation entre le \textit{Traité} et les travaux de Galois.\footnote{Pour un panorama général sur la relation Jordan-Galois voir : BRECHENMACHER (Frédéric), \href{http://smf4.emath.fr/Publications/RevueHistoireMath/17/html/}{Self-portraits with Évariste Galois (and the shadow of Camille Jordan)}, \textit{Revue d'histoire des mathématiques}, t. 17, fasc. 2, p. 271-369.} \textcolor{red}{Ainsi, Émile Picard déclarait-il à l'Académie des sciences le 23 janvier 1922 :}
\begin{quote}
\textcolor{red}{C'est surtout dans la théorie des substitutions et des équations algébriques que Jordan laisse une trace profonde. Dans un Ouvrage considérable sur les \textit{Substitutions}, il a fait une étude approfondie des idées de Galois [...]. Ces études ont permis à Jordan de résoudre un problème posé par Abel, celui de rechercher les équations de degré donné résolubles par radicaux et de reconnaître si une équation rentre ou non dans cette classe.\footnote{PICARD (Émile), Résumé des travaux mathématiques de Jordan, \textit{Comptes rendus de l'Académie des sciences de Paris}, t. 174 (1922), p. 210-211.}}
\end{quote}

Plus généralement, les travaux de Jordan ont souvent été placés dans le cadre d'une évolution de l'algèbre envisagée comme passant d'une théorie des équations à une discipline centrée sur des structures abstraites comme les groupes. Successeur de Jordan au Collège de France, Henri Lebesgue a publié une \textcolor{red}{biographie de notre héros qui mettait elle aussi en avant la relation Jordan-Galois} :
\begin{quote}
\textcolor{red}{Dans ses recherches, Jordan utilise la géniale méthode de Galois, dont le point essentiel est l'introduction d'un certain nombre de substitutions, déjà aperçu par Lagrange, que l'on peut attacher à chaque équation algébrique et dans lequel les propriétés des équations se reflètent fidèlement. Mais pour savoir observer dans ce miroir, il faut avoir appris à distinguer les diverses qualités des groupes de substitutions et à raisonner sur elles. C'est ce qu'a fait Jordan avec une habile ténacité et un rare bonheur ;  dans son \textit{Traité des Substitutions et des Équations algébriques}, où il a réuni et coordonné ses recherches, les propriétés des équations dérivent tout de suite de celles des groupes de substitutions. [...] Le théorème de Jordan sur la composition des groupes est le plus connu de tous ses résultats ; il entraîne cette conséquence fondamentale : il n'y a pas lieu de choisir entre les différents procédés de résolution algébrique d'une équations : ils sont tous équivalents et conduisent aux mêmes calculs à, l'ordre près.}\footnote{DIEUDONNÉ (Jean),  Notes sur les travaux de Camille Jordan relatifs à l'algèbre linéaire et multilinéaire et la théorie des nombres, \cite[p.V-XX]{9}.}
\end{quote}

Mais cette relation directe entre les idées de Galois et le \textit{Traité} de Jordan amène une tension entre les dimensions individuelles et collectives des travaux de ce dernier.
D'une part, Jordan a été célébré comme l'un des principaux fondateurs de la théorie des groupes, c'est à dire d'une théorie qui allait jouait un rôle majeur et structurant pour les organisation collectives des mathématiques au XX\up{e} siècle. Dans le cadre de l'édition des \oe{}uvres de Jordan en 1961, Jean Dieudonné présentait ainsi notre héros comme un second père de la théorie des groupes : \textcolor{red}{après Galois, le géniteur, Jordan serait le guide qui fait gagner l'âge adulte} :
\begin{quote}
\textcolor{red}{La théorie des groupes finis a été le sujet de prédilection de Jordan [...]. Son \oe{}uvre dans ce domaine est immense par le volume comme par l'importance, et son influence sur les développements ultérieurs de la théorie ne peut guère se comparer qu'à celle  des travaux de Galois lui-même. [...] Au moment où Jordan commence à écrire, la théorie des groupes est encore dans l'enfance, et en fait ce n'est guère qu'avec la publication de son traité qu'elle accèdera au rang de discipline autonome.}\cite{9}
\end{quote}

Mais d'autre part, et pour la même raison, Jordan a été présenté comme isolé et incompris en son temps. \textcolor{red}{Ainsi, dans la préface du premier tome des \OE uvres de Jordan, Gaston Julia estimait-il que :}
\begin{quote}
\textcolor{red}{Longtemps, Jordan a travaillé dans une solitude presque totale. Rares étaient ceux qui pouvaient apprécier la valeur de son \oe{}uvre. Aujourd'hui ses travaux sont plus actuels que lorsqu'ils ont été écrits, on les voit dans leur vraie lumière avec leur véritable portée. Dans cette lumière, Jordan nous apparaît, avec Galois et Sophus Lie, comme un des trois grands créateurs de la théorie générale des groupes. Cette reconnaissance de la valeur de ses méthodes et de ses résultats l'aurait certes réjoui, mais probablement pas surpris, car s'il était très modeste, il savait à coup sûr que son \oe{}uvre était solide et qu'elle porterait ses fruits dans l'avenir.}\cite{9}
\end{quote}

Les commentaires sur le \textit{Traité} présentent souvent une tension analogue entre l'individuel et le collectif. Cet ouvrage était ainsi présenté par le mathématicien américain James Pierpont à la fois comme une "synthèse unifiant les résultats des prédécesseurs de Jordan" et comme "contenant une somme immense de nouveaux résultats".\footnote{PIERPONT (James), The History of Mathematics in the Nineteenth Century, \textit{Bulletin of the American Mathematical Society}, vol. 11 , n° 2, 1904, p. 136-159.} Pour Bartel van der Waerden, la clef de voûte de cet alliage entre synthèse de travaux collectifs et création individuelle tiendrait à la structure même de l'ouvrage, présentée comme un "chef d'\oe{}uvre d'architecture mathématique, un édifice d'une beauté admirable".\footnote{WAERDEN (Bartel, van der), \textit{A history of Algebra : from Al-Khwàrizmi to Emmy Noether}, New York, Springer Verlag, 1985} Les commentateurs ne nous renseignent cependant guère sur cette structure. Au contraire, Félix Klein, qui se trouvait avec Sophus Lie à Paris au moment de la parution du \textit{Traité}, employait une métaphore biblique en comparant le \textit{Traité} à un "livre des sept sceaux".\footnote{KLEIN (Felix), \textit{Gesammelte mathematische Abhandlungen}, Springer: Berlin, 1921.}. 

L'inscription exclusive du \textit{Traité} dans l'émergence de la théorie de Galois et de la théorie des groupes donne donc une aura mystérieuse aux dimensions collectives des travaux de Jordan. Nous allons tenter dans cet article de dissiper ce mystère en analysant la structure du \textit{Traité}. 

Il faut tout d'abord être attentif à la dimension publique des commentaires sur la relation Jordan-Galois. De fait, ce n'est qu'au tournant des XIX\up{e} et XX\up{e} siècles que de telles appréciations apparaissent, c'est à dire à l'époque où la figure de Galois a été érigée en icône mathématique comme l'a étudié Caroline Ehrhardt.\footnote{EHRHARDT (Caroline), \textit{Evariste Galois. La fabrication d'une icône des mathématiques}, Paris, Éditions de l'EHESS, 2011.}  C'était en particulier à l'occasion du centenaire de l'École normale supérieure et de la réédition des \OE{}uvres de Galois que Sophus Lie et Émile Picard avaient attribué à Jordan un rôle majeur pour la mise en lumière des idées de Galois.\footnote{LIE (Sophus), Influence de Galois sur le développement des mathématiques, in DUPUY (Paul) (ed.), \textit{Le Centenaire de l'École Normale 1795-1895}, Paris, Hachette, 1895.\\
GALOIS (Évariste), \textit{\OE{}uvres mathématiques d'Évariste Galois, publiées sous les auspices de la Société Mathématique de France, avec une introduction par M. Émile Picard}, Paris: Gauthier-Villars, 1897.}

Un article disponible sur ce site (lien) précise le rôle historique attribué à Jordan en regard des rôles d'autorités endossés par des mathématiciens comme Picard.\footnote{Pour une discussion plus détaillée de ces questions, voir : BRECHENMACHER (Frédéric), "\href{ http://hal.archives-ouvertes.fr/aut/Frederic+Brechenmacher/}{Galois Got his Gun}," prépublication.} Ces discours impliquent souvent des hiérarchisations entre praticiens des mathématiques, notamment entre chercheurs et enseignants ou enseignants et ingénieurs. Les principales catégories mobilisées ("algèbre", "analyse", "théorie des équations", "théorie des groupes", "France", "Allemagne") sont ainsi loin d'être neutres. Elles participent au contraire de l'établissement de frontières entre mathématiques pures et appliquées, entre mathématiques élémentaires et supérieures, entre point de vue unificateur et réducteur, etc. Aux alentours de la Première Guerre mondiale, ces lignes épousent souvent celles des frontières nationales. Galois était alors souvent présenté comme un symbole  de l'universalité de la "pensée française". Un tel entrelacement de valeurs mathématiques et patriotiques se manifeste particulièrement \textcolor{red}{dans la nécrologie de Jordan écrite par Robert d'Adhémar en 1922 :}
\begin{quote}
\textcolor{red}{Jordan s'applique, dès 1860, à l'Algèbre de l'ordre, l'Algèbre des idées, bien plus haute que l'Algèbre des calculs, et, tout naturellement, il continue l'\oe{}uvre de cet enfant génial et décevant, Galois, qui, blessé dans un duel ridicule, mourut en 1832, âgé de 21 ans.
En 10 ans, Jordan construit ce qu'on appelle les Groupes des équations résolubles par radicaux et classe les équations non résolubles, distinguant celles qu'on peut ramener à des équations auxiliaires. Ses découvertes ont été publiées en 1870, dans le \textit{Traité des substitutions et des équations algébriques}, qui marque, après Abel et Galois, un progrès immense de l'Algèbre. [...] Chaque fois qu'il manie un être mathématique, Jordan met sur lui sa griffe puissant et austère. Là où il a été, la tranchée a été nettoyée ! }\footnote{ADHEMARD (Robert d'), Nécrologie. Camille Jordan, \textit{Revue générale des sciences pures et appliquées}, t. 3 (1922), p. 65-66.}
\end{quote}

Envisagée sur le temps long du XX\up{e} siècle, la structure des histoires présentées par de tels discours manifeste une certaine stabilité. Toutes présentent en effet la relation exclusive de Jordan aux idées de Galois comme donnant aux idées de ce dernier une dimension collective. Jordan était ainsi célébré pour avoir fait passer les travaux de Galois de l'état d'"idées" à une dimension "publique", rendant ces travaux disponibles "au monde entier".\footnote{PIERPONT (James), Early History of Galois Theory of Equations, \textit{Bulletin of the American Mathematical Society}, vol. 2, n° 4 (1897), p. 332-340.} La nature même de ces dimensions collectives change néanmoins radicalement entre les années 1900-1930 et 1930-1960, en passant de l'"analyse française" à l'"algèbre allemande" (lien vers l'autre article).

\section{Un portrait kaléidoscopique du \textit{Traité} dans les textes mathématiques entre 1870 et 1914}

À la différence des discours publics, les textes mathématiques publiés sur la période 1870-1914 font très rarement référence au \textit{Traité} en relation avec les travaux de Galois sur les équations algébriques. De fait, l'approche spécifique de Jordan sur la théorie de Galois n'a en réalité pas été reprise, au contraire d'autres approches, comme celle de  Joseph-Alfred Serret dans son \textit{Cours d'algèbre supérieure} de 1866 ou l'approche de Leopold Kronecker, radicalement opposée à celle de Jordan.\footnote{Voir à ce sujet : BRECHENMACHER (Frédéric), \href{http://smf4.emath.fr/Publications/RevueHistoireMath/17/html/}{La controverse de 1874 entre Camille Jordan et Leopold Kronecker}, \textit{Revue d'Histoire des Mathématiques},  13 (2007), p. 187-257.}  C'était notamment sur les travaux de Kronecker que s'appuyait Picard dans ses textes mathématiques et ce bien que les discours publics de ce dernier célébraient Jordan pour sa présentation des travaux de Galois.

De nombreux historiens ont interprété cette situation sous l'angle d'une faible réception du \textit{Traité} par les mathématiciens contemporains.\footnote{Voir notamment :\\
WUSSING (Hans), \textit{The Genesis of the Abstract Group Concept}, MIT Press, Cambridge: Mass.,1984.
KIERNAN (Melvin), The Development of Galois Theory from Lagrange to Artin, \textit{Archive for History of Exact Sciences}, vol. 8, n° 1-2 (1971), p. 40-152.}
Mais le \textit{Traité} est un ouvrage long et complexe. Il faut en envisager des réceptions multiples ainsi que des lectures fragmentées. 

L'impression d'isolement de Jordan provient principalement de perspectives rétrospectives ayant recherché les traces des travaux de ce dernier en théorie des groupes, théorie des équations, théorie de Galois ou plus généralement en algèbre. Or, non seulement ces catégories ne sont pas stables et revêtent des significations changeantes mais nous avons vu également que ces catégories sont loin d'être neutres et étaient au contraire employées par des autorités pour discourir publiquement des dimensions collectives des mathématiques. 

Comme nous le proposons en \hyperref[Annexe1]{\textcolor{blue}{annexe 1}}, dissocier la réception du Traité de la question du développement de la théorie des groupes ou de la théorie de Galois permet de brosser un panorama général de la réception des travaux de Jordan.\footnote{BRECHENMACHER (Frédéric), \href{http://hal.archives-ouvertes.fr/aut/Frederic+Brechenmacher/  }{On Jordan's measurements}, prépublication.} À la différence du portrait public de Jordan dans l'ombre de la figure de Galois, les textes mathématiques nous présentent un kaléidoscope de reflets mobiles et fragmentés qui nous donnent l'occasion d'éclairer les premiers travaux de Jordan sous un nouveau jour.

 Nous allons à présent porter notre attention sur un théorème mis en avant par certains lecteurs du \textit{Traité} et qui faisait déjà l'objet de la première thèse de Jordan : l'"origine" du groupe linéaire.

\section{Le premier théorème de Jordan : l'origine du groupe linéaire}

La première thèse présentée par Camille Jordan à la faculté des sciences de Paris en 1860 était consacrée au problème du "nombre des valeurs des fonctions". 

\subsection{Le problème du nombre des valeurs d'une fonction}

Ce problème s'avère l'une des racines de la théorie des groupes de substitutions. Il avait émergé de travaux du XVIII\up{e} siècle sur les équations. Comme nous l'illustrons en \hyperref[Annexe2]{\textcolor{blue}{annexe 2}}, la résolubilité par radicaux d'une équation algébrique de degré $n$ avait été mise en relation avec le nombre de valeurs qu'une fonction "résolvente" de $n$ variables peut prendre lorsque ses variables sont permutées de toutes les manières possibles. 

Etant donnée une fonction $\phi(x_1, x_2, ..., x_n)$ de $n$ "lettres", une "valeur" de $\phi$ est une fonction obtenue par permutation des variables pour toute substitution $\sigma\in\Sigma(n)$,
\[
\phi^\sigma(x_1, x_2, ..., x_n)=\phi(x_{1\sigma}, x_{2\sigma}, ..., x_{n\sigma})
\]

En général, une fonction peut ainsi prendre $n!$ valeurs différentes mais il peut arriver que quelques unes de ces valeurs deviennent identiques et le problème abordé par la thèse de Jordan consiste à déterminer le nombre de valeurs pouvant être prises pour certaines classes de fonctions. En termes contemporains, le problème revient à déterminer tous les ordres possibles des sous-groupes du groupe symétrique $\Sigma(n)$.

Au cours du XIX\up{e} siècle, ce problème avait donné lieu à des développements autonomes de la théorie des équations. Augustin-Louis Cauchy (1815, 1844-46), Joseph Bertrand, (1845) et Serret (1849) avaient notamment énoncé des bornes pour le nombre de valeurs pouvant être prises par des fonctions de certains types généraux.\footnote{Dès le début du XIX\up{e} siècle, Ruffini et Cauchy avait par exemple montré que  le nombre de valeurs d'une fonction non symétrique de cinq variables ne peut être moindre de cinq à moins qu'il soit égal à deux.}  Charles Hermite et Kronecker avaient davantage porté leur attention sur des fonctions spéciales, comme une fonction de 6 variables prenant exactement 6 valeurs. 
À partir des années 1850, ces travaux pouvaient s'appuyer sur un développement théorique dans lequel Cauchy avait mis en avant la notion essentielle de "système de substitutions conjuguées" (c'est-à-dire ce que nous appelons des groupes finis et leurs sous-groupes distingués). En 1860, le problème du nombre de valeurs des fonctions avait été choisi comme sujet du Grand prix des sciences de mathématiques de l'Académie des sciences.\footnote{Voir au sujet de cet épisode, ainsi que des travaux de Thomas Kirkman, également candidat au grand prix de 1860,  EHRHARDT (Caroline), \textit{Evariste Galois et la théorie des groupes. Fortune et réélaborations (1811-1910)}, Thèse de doctorat. Ecole des Hautes études en sciences sociales. Paris, 2007, p. 291-393.}  Deux jeunes mathématiciens y consacraient leurs premières recherches : Jordan et Émile Mathieu.

\subsection{Une approche générale par des réductions successives}

Dans sa thèse, Jordan revendiquait une approche "générale" du problème par des "réductions successives" en sous-problèmes. La première réduction montrait que l'étude des fonctions générales transitives (correspondant à des équations algébriques irréductibles) peut se réduire à l'étude de fonctions sur des \textcolor{blue}{"systèmes imprimitifs"}\footnote{ \textcolor{blue}{D'un point de vue contemporain, il s'agit ici de décomposer les éléments du corps en blocs d'imprimitivité sous l'action d'un groupe de substitution imprimitif, lui-même décomposé en un groupe primitif quotient. Soit $G$ un groupe transitif opérant sur ensemble $\Omega$. Un sous-ensemble $\Gamma$ de $\Omega$ est appelé bloc d'imprimitivité si $\Gamma \ne \emptyset$ et si pour tout $g\in G$, alors soit $\Gamma g=\Gamma$ soit $\Gamma g \cap \Gamma = \emptyset$. Si $\Gamma$ est un tel bloc et $\Gamma_1,\Gamma_2, ..., \Gamma_m$,  sont les ensembles distincts obtenus par une orbite $\Gamma g$ pour $g \in G$, alors $(\Gamma_1, \Gamma_2, ..., \Gamma_m)$   est une partition de $\Omega$. $G$ est dit imprimitif s'il existe un bloc propre non trivial dans $\Omega$. $G$ est dit primitif s'il n'est pas imprimitif.}} de lettres : l'ensemble de toutes les lettres peut dans ce cas se diviser en \textcolor{blue}{blocs}\footnote{ \textcolor{blue}{Jordan désignait en réalité par le terme "groupe" les blocs de lettres, c'est à dire un sous-ensemble d'un corps fini,  et ce que nous appellerions des groupes des substitutions par le terme de système conjugué introduit par Cauchy. Rappelons que chez des auteurs antérieurs, comme Galois ou Poinsot, le terme groupe concernait intrinsèquement à la fois des arrangements de lettres en blocs de permutations et les opérations sur de tels arrangements par des substitutions. Voir notamment à ce sujet :
DAHAN (Amy), Les travaux de Cauchy sur les substitutions. Étude de son approche du concept de groupe, \textit{Archive for History of Exact Sciences}, vol. 23, 1980, p. 279-319.}} - représentés ci-dessous par une succession de lignes - de manière à ce que les substitutions opèrent soit en permutant entre elles les lettres d'une même ligne soit les lignes les unes avec les autres :
\[
\begin{matrix}
a_1 & a_2 & ... & a_m\\
b_1& b_2 & ... & b_m\\
c_1 & c_2 & ... & c_m\\
. & . & . & .
\end{matrix}
\]

L'enjeu de cette première réduction était d'introduire une indexation des lettres par des suites d'entiers $(1, 2, ..., m)$ de manière à faire opérer les substitutions non plus sur des lettres mais sur des nombres. 

\subsection{La représentation analytique des substitutions}

Une telle indexation était un préalable à l'introduction du problème crucial de la représentation analytique des substitutions : étant donnée une  substitution $S$ sur $p$ lettres $a_0$, $a_1$, ..., $a_{p-1}$, ce problème consistait à trouver une fonction analytique $\phi$ telle  que 
\[
S(a_i)=a_{\phi(i)}
\]
Bien qu'elle soit passée inaperçue de nombreux travaux historiques, la représentation analytique des substitutions a joué un rôle important au XIX\up{e} siècle.

Comme le détaille l'\hyperref[Annexe3]{\textcolor{blue}{annexe 3}}, cette représentation avait notamment été mise en avant au début du XIX\up{e} siècle par  Louis Poinsot à la suite des travaux de Carl Friedrich Gauss sur les équations cyclotomiques.\footnote{Voir à ce sujet : BOUCARD (Jenny), \href{http://smf4.emath.fr/Publications/RevueHistoireMath/17/html/}{Louis Poinsot et la théorie de l'ordre : un chaînon manquant entre Gauss et Galois ?}, \textit{Revue d'histoire des mathématiques}, 17, fasc. 1 (2011), p. 41-138. FREI (Gunther), The unpublished section eight: On the way to function fields over a finite field, in GOLDSTEIN (Catherine), SCHAPPACHER (Norbert), SCHWERMER (Joaquim) (eds.), \textit{The Shaping of Arithmetics after C. F. Gauss's Disquisitiones Arithmeticae}, Berlin: Springer, 2007, p. 159-198.
NEUMANN (Olaf), The Disquisitiones Arithmeticae and the Theory of Equations, in GOLDSTEIN, SCHAPPACHER, SCHWERMER,  \textit{op. cit.}, p. 107-128.} Étant donné un nombre $p$ premier, l'équation cyclotomique associée à l'équation binôme 
\[
x^p-1=0
\]
est l'équation irréductible : 
\[
\frac{x^p-1}{x-1}=x^{p-1}+x^{p-2}+...+x+1
\]

Toutes les racines de cette équation peuvent être exprimées en fonctions de l'une d'entre elles, $\omega$ (une racine primitive) :
\[
\omega, \omega^2,..., \omega^{p-1}
 \]
 
Les substitutions des racines forment ainsi un groupe cyclique engendré par une substitution pouvant être représentée analytiquement par son action sur les puissances $i$ indexant les racines , c'est à dire par $(i \ i+1)$. 
Une propriété cruciale des équations cyclotomiques est que la suite des racines peut être réindexée par une racine primitive $g$ de l'équation binôme de congruence : 
\[
i^{p-1}-1 \equiv 0 \ mod(p)
\]
de manière à obtenir la suite :
\[
\omega^g, \omega^{g^2},..., \omega^{g{p-1}}
 \]
Le groupe cyclique est alors engendré par une substitution représentée analytiquement par $(i \ gi)$.

\subsection{Le premier théorème de Jordan}

Mais revenons à la première thèse de Jordan. Après avoir réduit le problème à l'étude des fonctions imprimitives, Jordan montrait comment poursuivre la réduction aux fonctions primitives, c'est à dire au cas où il n'est pas possible de subdiviser les lettres en plusieurs lignes, $\Gamma_1, \Gamma_2, ..., \Gamma_m$,  que les substitutions permuteraient en blocs. Le point clef de cette réduction consistait à montrer que les substitutions opérant sur les systèmes imprimitifs peuvent être décomposées en deux espèces, correspondant aux deux formes de représentations analytiques des cycles sur le modèle de la décomposition réalisée par Gauss pour les équations cyclotomiques (\hyperref[Annexe3]{\textcolor{blue}{annexe 3}}).

\begin{itemize}
\item  La première espèce permute cycliquement les blocs  $\Gamma_1, \Gamma_2, ..., \Gamma_m $ eux-mêmes par des substitutions du type $(i \ i+1)$ sur les indices. Dans le cas plus général de lettres indexées par $n$ indices $a_{x, x', x'',...}$ ces substitutions prennent la forme :
\[
a_{x+ \alpha \ mod.p, \ x'+\alpha' \ mod.p, \ x''+\alpha'' \ mod.p, \ ...}
\]
\item La seconde espèce permute cycliquement les lettres à l'intérieur de chaque bloc $\Gamma_i$ par des substitutions du type $(i \ gi)$. Dans le cas de $n$ indices, elles prennent ainsi la forme :

\[
a_{ax+bx'+cx''... mod. p, a'x+b'x'+c'x''... mod.p, a''x+b''x'+c''x''... mod.p}
\]
\end{itemize}

Jordan énonçait alors son \textcolor{blue}{principal théorème} :\footnote{\textcolor{blue}{En termes contemporains, le théorème de Jordan peut s'interpréter comme énonçant que le sous-groupe normal minimal $A$ d'un groupe primitif résoluble $G$ est abélien du type $(1, 1,..,1)$ (donc isomorphe à un produit de groupes cycliques ou au groupe multiplicatif d'un corps fini $GF_{p^n}*$). L'action de $G$ sur $A$ définit un groupe linéaire.}} 

\begin{quote}
\textbf{Premier théorème de Jordan}
\begin{itemize}
\item Les systèmes primitifs comportent un nombre de lettres donné par la puissance d'un nombre premier $p^n$.
\item Les substitutions sur ces systèmes ont une forme analytique linéaire : 
\[
(x, x', x'' ... ; ax+bx'+cx''+... +d, a'x+b'x'+c'x''+...+d', a''x+b''x'+c''x''+..+d'', ...)
\]
que l'on peut également dénoter par : 
\[
\begin{vmatrix}
x & ax+bx'+cx''+... \\
x' & a'x+b'x'+c'x''+... \\
x'' & a''x+b''x'+c''x''+... \\
.. & .....................
\end{vmatrix}
\]
\end{itemize}
\end{quote}
Ce groupe devait être désigné par Jordan quelques années plus tard sous le nom de \textcolor{blue}{"groupe linéaire"}.\footnote{Nous le désignerions aujourd'hui comme un groupe affine sur un corps fini.}. 

\subsection{Retour sur la relation Jordan-Galois}

C'est en déterminant l'ordre de ce groupe que Jordan s'apercevait de la relation entre ses propres travaux et ceux menés par Galois pour le cas de $p^2$ lettres (\hyperref[Annexe3]{\textcolor{blue}{annexe 3}}). Dans la version modifiée de sa thèse, publiée en 1861, Jordan ajoutait alors un supplément dans lequel il commentait le "Fragment de second mémoire" de Galois sur les \textcolor{blue}{équations primitives de degré $p^n$}.\footnote{ \textcolor{blue}{Une équation irréductible de degré composé $n=pq$ a un bloc d'imprimitivité de taille $p$ si et seulement si il existe une équation auxiliaire $g(x)=0$ de degré $m$ telle que l'adjonction de toutes les racines de cette équation permet de factoriser $f(x)$ en un produit de $p$ facteurs irréductibles de degré $q$. Considérons par exemple $f(x)=x^6-2=0$ et $g(x)=x^3-2$. Le degré du corps de décomposition $K$ de $f$ sur $\Q$ est 12 et le groupe de Galois $G$ est le groupe dihédral d'ordre 12, qui peut être représenté comme le groupe des symétries d'un hexagone régulier (dont les sommets sont indexés par les racines sixièmes de 2). $G$ est engendré par deux éléments : $x$, la rotation d'angle $\frac{\pi}{3}$, et $y$, la symétrie d'axe une diagonale. Soit alors les deux générateurs de $K$ : $\alpha=\sqrt[6]{2}$ et $\omega=e^{\frac{2 \pi}{3}}$, les éléments de $G$ peuvent être décrits comme automorphismes de $K$ par :
\[
 x: \alpha \rightarrow e^{\frac{\pi}{3}\alpha} ,  \omega \rightarrow \omega \ ; \ 
 y : \alpha \to \alpha, \omega \to \omega^2 
\]
Alors $G_\alpha$ le groupe cyclique d'ordre 2 engendré par$y$ (i.e. la conjugaison complexe) est un bloc d'imprimitivité. Si l'on adjoint à $\Q$ la racine $\gamma=\sqrt[3]{2}$ de $g(x)=0$ ainsi qu'une racine primitive $\omega$ de l'unité, on obtient en effet la factorisation : $f(x)=(x^2-\gamma)(x^2-\omega \gamma)(x^2-\gamma^2 \alpha)$ sur $\Q(\gamma, \omega)$. Voir à ce sujet : NEUMANN (Peter M.),The concept of Primitivity in Group Theory and the Second Memoir of Galois, Archive for History of Exact Sciences, 60 (2006), p. 379-429.}}

Cet intérêt de Jordan pour la question générale de la résolubilité des équations de degré composé, et donc des groupes primitifs de degré $p^n$, s'avère la principale spécificité de la relation établie par ce dernier aux travaux de Galois. Contrairement à d'autres aspects des travaux de Galois, comme le critère de résolubilité des équations de degré premier, les imaginaires de la théorie des nombres ou les équations modulaires, le cas des équations primitives de degré $p^n$ n'avait en effet avant Jordan été envisagé que par Enrico Betti et Alexandre Allégret (\hyperref[Annexe3]{\textcolor{blue}{annexe 3}}). 

Jordan s'appuyait ainsi principalement sur une lecture du "Fragment du second mémoire" de Galois et non sur le célèbre "Mémoire sur les conditions de résolubilité des équations par radicaux" dans lequel se trouvent les principes généraux de ce que nous appelons aujourd'hui la théorie de Galois. Au contraire de Jordan qui, comme nous l'avons vu, s'intéressait ainsi aux groupes linéaires de $n$ variables, la grande majorité des présentations des travaux de Galois se focalisaient jusqu'à la fin du siècle sur les principes généraux du "Mémoire" et leur application aux équations de degré premier $p$ qui nécessitent uniquement de considérer des groupes linéaires à une variable $(i \ ai+b)$.

Comme nous l'avons vu, Jordan avait cependant introduit le groupe linéaire dans sa thèse indépendamment de Galois et ce n'était qu'a posteriori qu'il s'était intéressé aux travaux de ce dernier sur la résolubilité des équations. Dans ses travaux ultérieurs, Jordan insistait à de nombreuses reprises sur la différence entre ses propres travaux et ceux de Galois. Il critiquait notamment la généralisation abusive par Galois de son critère de résolubilité des équations primitives de degré $p$ aux équations primitives de degré $p^n$ (\hyperref[Annexe3]{\textcolor{blue}{annexe 3}}). \textcolor{red}{Galois avait en effet énoncé que de telles équations sont résolubles si et seulement si leur groupe est linéaire} :
\begin{quote}
\textcolor{red}{Galois avait annoncé que les équations primitives et solubles par radicaux rentreraient dans un type unique, sauf pour le neuvième et le vingt-cinquième degré, qui présenteraient certains types exceptionnels. On voit par les énoncés qui précèdent qu'il faut prendre presque exactement le contre-pied de cette assertion.}\cite[p.113]{3}
\end{quote}

Or la forme linéaire des substitutions n'est pas un critère suffisant pour caractériser les groupes résolubles de degré non premier. Comme le montrait Jordan au cours des années 1860, il est nécessaire de poursuivre plus avant la chaîne successive de réductions du problème en recherchant les sous-groupes résolubles du groupe général linéaire.\\

Récapitulons. Comme nous l'avons vu, l'objet principal de la thèse de Jordan était d'énoncer un théorème "engendrant" le groupe linéaire. Ce groupe n'était ainsi pas défini, comme nous le ferions aujourd'hui, par une liste d'axiomes. Au contraire, le groupe linéaire se présentait comme engendré par une chaîne de réductions successives d'un problème général en sous-problèmes.

Par ailleurs, les substitutions linéaires étaient avant tout identifiées par leur forme analytique. Comme nous le détaillons en \hyperref[Annexe3]{\textcolor{blue}{annexe 3}} et comme nous allons le voir plus précisément dans la suite, la représentation analytique des substitutions n'était pas une simple notation mais s'accompagnait de procédures spécifiques de réductions auxquelles Jordan attribuait l'"essence" d'une approche transversale à diverses théories mathématiques.\\

\section{L'"essence" du Traité des substitutions et des équations algébriques}

Abordons à présent la question de la structure du \textit{Traité} de 1870. Nous allons retrouver le premier théorème de Jordan dispersé dans les quatre livres qui composent cet ouvrage. Mais nous allons voir également que ce théorème donne l'"essence" sous-jacente à la structure du \textit{Traité}. 

La réduction de la représentation analytique des substitutions linéaires par les deux formes de représentations de cycles sous-tend en effet une chaîne de généralisations successives courant les trois premiers livres du \textit{Traité} jusqu'à l'énoncé du "théorème fondamental" sur la résolubilité des équations algébriques. Ce théorème permet alors de retourner cette chaîne en une suite de réductions successives d'un problème général.\\

\subsection{Livre I. Des congruences}

L'ouverture du \textit{Traité} introduit les notions relatives aux congruences qui permettent l'indexation de lettres et, par là,  la représentation analytique des substitutions. Après l'exposé des propriétés des congruences binômes, ce premier livre expose ce que Jordan dénommait la "théorie de Galois". Cette désignation ne correspond pas à ce que recouvre aujourd'hui cette théorie mais désigne les  "imaginaires  de Galois" permettant d'indexer des systèmes de $p^n$ lettres en généralisant les propriétés des racines cyclotomiques des équations de degré $p$ aux congruences irréductibles de degré $n \ mod.p$  (\hyperref[Annexe3]{\textcolor{blue}{annexe 3}}). 
Ce premier livre est ainsi l'occasion d'introduire un cas particulier de représentation analytique de substitutions : les cycles $(i \ i+a)$ et $(i \ gi)$.

\subsection{Livre II. Des substitutions}
Un premier chapitre de généralités sur les substitutions présente une synthèse de résultats antérieurs de mathématiciens comme Cauchy, Serret,  Bertrand ou  Mathieu. Le second chapitre, qui représente la partie la plus importante de ce deuxième livre, est consacré aux propriétés des groupes linéaires généraux et spéciaux étudiés par Jordan depuis sa thèse. 

Le problème de la représentation analytique des substitutions ouvre ce second chapitre et amène la "génération du groupe linéaire". Nous avons vu que l'origine de ce groupe avait fait l'objet du premier théorème de Jordan. Le théorème se présente cependant la tête en bas en 1870 : alors que le groupe linéaire était engendré en 1860 par des réductions successives du problème du nombre de valeurs d'une fonction, cette génération intervient à présent à partir du problème de la recherche de la forme analytique de l'ensemble des substitutions $g$ \textcolor{blue}{"permutant" des produits de cycles $(x \ x+1)$ de la forme analytique suivante} :\footnote{\textcolor{blue}{Il s'agit, en termes contemporains, d'un groupe abélien élémentaire.}}
\[
 |x, x', ..., x+\alpha, x'+\alpha', ... |
\]
\textcolor{blue}{c'est à dire transformant une substitution élémentaire $c$ de la forme ci-dessus en une autre $c'$ : $gcg-1=c'$} :\footnote{\textcolor{blue}{En termes contemporains, il s'agit de rechercher le plus grand groupe dans lequel un groupe abélien élémentaire est distingué.}}
\begin{quote}
Ces dernières substitutions, que nous représenterons indifféremment par l'une ou l'autre des deux notations suivantes :
\[
\begin{vmatrix}
x & ax+bx'+cx''+... \\
x' & a'x+b'x'+c'x''+... \\
x'' & a''x+b''x'+c''x''+... \\
.. & .....................
\end{vmatrix}
, (x, x', ... ; ax+bx'+ ..., a'x+b'x'+ ..., ...)
\]
forment évidemment un groupe, que nous appellerons le groupe linéaire. \cite[p.91]{6}
\end{quote}

La démonstration, détaillée ici (lien), consiste à générer le groupe linéaire par les deux formes de représentations analytiques des cycles introduites au Livre I.

\subsection{Livre III. Des irrationnelles}

Le troisième livre s'ouvre par une présentation des principes généraux de Galois sur la résolubilité des équations algébriques. Ce premier chapitre du Livre III a été amplement commenté, je ne vais pas le détailler davantage ici.\footnote{Voir à ce sujet : EHRHARDT (Caroline), \textit{Evariste Galois et la théorie des groupes. Fortune et réélaborations (1811-1910)}, Thèse de doctorat. Ecole des Hautes études en sciences sociales. Paris, 2007.
CORRY (Leo), \textit{Modern Algebra and the Rise of Mathematical Structures}, Basel: Birkhäuser, 1996.
KIERNAN (Melvin), The Development of Galois Theory from Lagrange to Artin, \textit{Archive for History of Exact Sciences}, vol. 8, n° 1-2 (1971), p. 40-152}  Notons l'expression "Des irrationnelles" utilisée par Jordan pour ce que nous désignons aujourd'hui comme la théorie de Galois et sur laquelle  \hyperref[Annexe4]{\textcolor{blue}{l'annexe 4}} propose un éclairage.

Le livre se poursuit par trois chapitres d'applications : applications algébriques, applications géométriques et applications aux transcendantes. Le chapitre consacré aux équations algébriques éclaire sous un nouveau jour l'origine du groupe linéaire telle que présentée au Livre II comme permutant les cycles introduits au Livre I. 

Jordan commence par considérer le cas des équations abéliennes, dont toutes les racines s'expriment en fonction rationnelle de l'une d'entre elles, et plus spécifiquement le cas des équations cyclotomiques pour lesquelles les racines sont donnés par la suite des puissances d'une raine primitive. Dans ce cas, comme nous l'avons vu, le groupe associé à l'équation est engendré par le cycle $(i \ i+1)$ mais les racines peuvent être réindexées de manière ce que le groupe cyclique soit engendré par le cycle $(i \ gi)$.
Le cas des "équations de Galois", pour lesquelles toutes les racines sont fonctions rationnelles de deux d'entre elles, se présente alors comme une généralisation naturelle des équations abéliennes sur le modèle de l'origine du groupe linéaire. En effet, ces équations sont associées à un groupe de substitutions linéaires de la forme $(i \ ai+b)$ (\hyperref[Annexe3]{\textcolor{blue}{annexe 3}}). Or le groupe linéaire a précisément été engendré au Livre II par la combinaison des cycles $(i \ i+1)$ et $(i \ gi)$. L'origine du groupe linéaire donne ainsi un modèle de réduction du groupe d'une équation de Galois: $(i \ ai+b)$ se décompose en effet en $(i \ i+1)$ et $(i \ gi)$. La résolution d'une équation de Galois se ramène par conséquent à celle de deux équations abéliennes. 

Comme nous allons le voir à présent, cette situation donne à son tour un modèle pour un ultime maillon de la chaîne de généralisation présentée par Jordan.

\subsection{Livre IV. De la résolution par radicaux}

Ce livre, qui représente à lui seul plus du tiers de l'ensemble de l'ouvrage, s'ouvre comme les précédents par un premier chapitre de généralités. Jordan y énonce une suite de théorèmes se concluant par le \textcolor{blue}{"théorème fondamental"}\footnote{\textcolor{blue}{En des termes qui nous sont contemporains, ce théorème énonce qu'une équation irréductible est résoluble par radicaux  si et seulement si son groupe de Galois peut être "réduit" en une suite $G=H_0 \supset H_1 \supset H_2 \supset ...\supset H_m = I$ dans laquelle chaque $H_k$ est un sous-groupe normal de $G$ et les quotients $\frac{H_k}{H_{k+1}}$ sont abéliens.}} suivant :
\begin{quote}
\textbf{Théorème fondamental du Livre IV} \\
Pour qu'un groupe $L$ soit résoluble, il faut et il suffit qu'on puisse former une suite de groupes $1, F, G, H, ..., L$ se terminant par $L$, et jouissant des propriétés suivantes : 1° chacun de ces groupes est contenu dans le suivant, et permutable aux substitutions de $L$; 2° deux quelconques de ses substitutions sont échangeables entre elles, aux substitutions près du groupe précédent \cite[p.395]{6}.
\end{quote}
Ce théorème se présente comme une généralisation de la décomposition du groupe des équations de Galois en deux groupes cycliques ; il englobe ainsi le critère de résolubilité énoncé par Galois pour les équations de degré premier (\hyperref[Annexe3]{\textcolor{blue}{annexe 3}}). 
Mais nous retrouvons surtout dans le théorème fondamental une dimension cruciale du premier théorème de Jordan, à savoir l'approche générale que ce dernier entendait développer au moyen d'une chaîne de réductions successives. 

L'énoncé de ce théorème permet alors de retourner la chaîne de généralisations que nous avons vu courir au long des trois premiers livres du \textit{Traité}. Le théorème fondamental sous-tend en effet une chaîne de réductions des groupes les plus généraux aux groupes les plus spéciaux - les groupes cycliques - en passant par les groupes linéaires. C'est précisément à une telle réduction que s'emploie le Livre IV en abordant le problème de la recherche des groupes résolubles généraux par des \textcolor{red}{réductions successives à des groupes résolubles spéciaux : groupes transitifs, primitifs, linéaires, symplectiques etc.}
\begin{quote}
\textcolor{red}{Supposons que nous ayons formé, pour un degré donné, le tableau de tous les groupes résolubles et transitifs les plus généraux. Chacun d'eux caractérisera un type distinct d'équations irréductibles résolubles par radicaux. Les groupes résolubles et transitifs, non généraux, caractériseront des types d'équations plus spéciaux, et contenus dans les précédents comme cas particuliers.}\cite[p. 396]{6}
\end{quote}
Plus encore, cette chaîne de réductions constituait, pour Jordan, l'essence même de sa méthode :
\begin{quote}
L'essence de ma méthode consiste à déterminer successivement les groupes partiels $F, G, H, ...$.\cite[p.963]{2}.
\end{quote}
Dans cette chaîne, le groupe linéaire joue un rôle essentiel en tant que groupe le plus général dont les substitutions disposent d'une représentation analytique. Le Livre IV fait en effet un usage constant des procédés de décompositions des formes analytiques des substitutions linéaires. 

\section{La théorie de l'ordre}

Nous avons vu le rôle joué dans le \textit{Traité} par une chaîne de réductions des groupes prenant modèle sur la réduction de la représentation analytique des substitutions linéaires en produits de cycles. Nous avons vu également que cette approche n'était pas directement liée à ce que nous appelons aujourd'hui "théorie de Galois" : non seulement Jordan avait développé de telles méthodes dès sa thèse, soit avant d'étudier les écrits de Galois, mais il avait à plusieurs reprises insisté sur la différence entre son approche et les principes généraux énoncés par Galois pour la résolubilité des équations. Par conséquent, les travaux de Jordan ne peuvent pas s'envisager uniquement sous l'angle d'une relation exclusive aux idées de Galois. 

Abordons à présent la question cruciale des dimensions collectives des recherches de Jordan dans les années 1860. Dès l'introduction de sa thèse, Jordan avait placé le problème du nombre des valeurs d'une fonction dans le cadre de la "théorie de l'ordre" en revendiquant l'héritage des travaux de Poinsot. \textcolor{red}{Cette théorie était présentée comme englobant les travaux de Gauss sur les équations cyclotomiques, ceux d'Abel et de Galois sur la résolubilité des équations ainsi que l'approche de Cauchy sur les déterminants :}
\begin{quote}
\textcolor{red}{Mais il pourra arriver que quelques-unes d'entre elles [les fonctions obtenues par permutation de leurs variables] deviennent identiques, par suite de quelque symétrie que présente la fonction primitive. L'étude de ces diverses sortes de symétrie offre un grand intérêt ; car c'est la base et le point de départ naturel de ce genre de recherches que M. Poinsot a distingué de tout le reste des mathématiques sous le nom de théorie de l'ordre : elle présente en outre d'importantes applications. C'est dans le Mémoire où M. Cauchy a donné les premiers principes généraux de cette théorie, qu'il a établi pour la première fois les théorèmes fondamentaux sur les déterminants. Abel s'est appuyé sur elle pour établir l'irrésolubilité de l'équation générale du cinquième degré. Galois, dans un admirable Mémoire, en a fait dépendre, non seulement les conditions de la résolution algébrique, mais la théorie entière des équations, considérée sous son point de vue le plus général, et la classification des irrationnelles algébriques.}\cite{1}
\end{quote}
Des travaux récents de Jenny Boucard ont montré que la théorie de l'ordre avait été envisagée par Poinsot comme articulant les analogies présentées par les situations cycliques rencontrées en algèbre, théorie des nombres, géométrie et mécanique.\footnote{BOUCARD (Jenny), \href{http://smf4.emath.fr/Publications/RevueHistoireMath/17/html/}{Louis Poinsot et la théorie de l'ordre : un chaînon manquant entre Gauss et Galois ?}, \textit{Revue d'histoire des mathématiques}, 17, fasc. 1 (2011), p. 41-138.}

Pour Jordan, la théorie de l'ordre était associée au caractère "essentiel" de sa "méthode de réduction" de la représentation analytique des substitutions. Cette "réduction" exprimait en effet le type de relations cachées entre diverses théories ou classes d'objets auxquelles Jordan attribuait l'essence de son approche. Ainsi, tout en admettant que sa réduction des groupes n'était pas des plus avantageuses dans les applications, Jordan insistait que "si l'on se borne à étudier le problème de la symétrie en lui-même, cette méthode, plus naturelle et plus directe, peut seule conduire aux véritables principes". \textcolor{red}{Et d'ajouter :}
\begin{quote}
\textcolor{red}{On pourrait voir une image de ce résultat dans le théorème de mécanique qui ramène le mouvement général d'un corps solide à un mouvement de translation combiné avec une rotation autour du centre de gravité. Ce principe du classement des lettres en divers groupes est le même dont Gauss et Abel ont déjà montré la fécondité dans la théorie des équations: il me semble être dans l'essence même de la question, et sert de fondement à toute mon analyse.}\cite[p.5]{1}
\end{quote}
Le procédé de réduction des substitutions linéaires en produits de cycles représentés par $(i \ i+1)$ et $(i \ gi)$ et sa généralisation par Jordan à la décomposition des groupes imprimitifs étaient ainsi envisagé comme une sorte de dévissage par analogie avec la décomposition du mouvement hélicoïdal d'un solide en mouvements de rotation et de translation. 
Jordan allait par la suite développer cette analogie dans ses travaux en s'inspirant des travaux géométriques et mécaniques de Poinsot mais aussi de ceux de Bravais sur la cristallographie.\footnote{SCHOLZ (Erhard), \textit{Symmetrie, Gruppe, Dualität. Zur Beziehung zwischen theoretischer Mathematik und Anwendungen in Kristallographie und Baustatik des 19. Jahrhunderts}, Basel: Birkhäuser, 1989.} Chez Jordan,  la théorie de l'ordre se présentait comme transversale à différents domaines qui allaient représenter l'essentiel des préoccupations de notre héros jusqu'en 1867-1868 :
\begin{itemize}
\item théorie des nombres (cyclotomie, congruences)
\item algèbre (équations, substitutions) 
\item analyse (groupes de monodromie et lacets d'intégration des équations différentielles linéaires)
\item géométrie/ topologie (cristallographie, symétries des polyèdres et des surfaces (y compris de Riemann))
\item mécanique (mouvements des solides).
\end{itemize}
 Lors de la présentation de ses travaux à l'occasion de sa candidature à l'Académie en 1881, Jordan décrivait l'orientation générale de ses travaux comme portant davantage sur les relations entre des classes d'objets que sur les objets eux-mêmes \textcolor{red}{dans l'héritage des travaux de Poinsot :}
 \begin{quote}
 \textcolor{red}{Les Mathématiques ne sont pas seulement la science des rapports, je veux dire que l'esprit n'y a pas seulement en vue la proportion et la \textit{mesure}; il peut encore considérer le \textit{nombre} en lui-même, l'\textit{ordre} et la \textit{situation} des choses, sans aucune idée de leurs rapports ni des distances plus ou moins grandes qui les séparent. Si l'on parcourt les différentes parties des Mathématiques, on y trouve partout ces deux objets de nos spéculations.
Ainsi, à côté de l'Algèbre ordinaire, il y a une Algèbre supérieure, qui repose tout entière sur la théorie de l'ordre et des combinaisons. D'autre part, ce qui rend la théorie des polyèdres très difficile, c'est qu'elle tient essentiellement à une science presque encore neuve, que l'on peut nommer \textit{Géométrie de situation}, parce qu'elle a principalement pour objet, non la grandeur ou la proportion des figures, mais l'ordre ou la situation des éléments qui les composent.
[...] la tendance générale de mes recherches [a] eu presque constamment pour but d'approfondir la \textit{théorie de l'ordre} au double point de vue de la Géométrie pure et de l'Analyse.
En Géométrie, j'ai étudié successivement les lois de la symétrie des polyèdres, des systèmes de lignes et des systèmes de molécules.
En Analyse, j'ai pris pour objet principal de mes travaux la théorie des substitutions (qui n'est au fond autre chose que celle de la symétrie des expressions algébriques) et ses applications à la théorie des équations algébriques et celle des équations différentielles linéaires.}\cite[p.7-8]{8}
 \end{quote}
Cette approche transversale permet de jeter un nouvel éclairage sur la dimension collective des travaux de Jordan des années 1860. En effet, ces travaux manifestent des interactions avec des publications de nombreux autres mathématiciens contemporains dans divers domaines comme les mouvements des solides et l'étude des polyèdres,\footnote{Bertini, Godt,  Goursat,  Kirkmann, Fedorow,  Hagen, etc.} la géométrie cinématique,\footnote{Houël, Méray, etc.}  la déformation des surfaces,\footnote{Becker, Boussinesq}  la cristallographie,\footnote{Sohncke}  les substitutions,\footnote{Allégret, Despeyroux, Mathieu.} ou encore l'étude des fonctions spéciales en lien avec les surfaces de Riemann et leurs applications aux équations différentielles.\footnote{Schwarz, Klein, Gierster, Hess, Hurwitz, Bukhardt etc.} 

Parmi ces mathématiciens, les auteurs français étaient pour la plupart liés à l'École polytechnique. On peut donc soupçonner la transmission d'une culture mathématique spécifique dans le cadre de l'enseignement polytechnicien de cette époque. Mais un domaine semblable à ce que Jordan désignait comme la théorie de l'ordre était également reconnu par d'autres acteurs à une échelle européenne. Ce domaine semble donc constituer un champ de recherche spécifique dans les années 1860-1880. La question des dimensions collectives d'un tel ensemble de travaux reste cependant encore ouverte et appelle davantage de recherches historiques.

L'identification de l'approche transversale de Jordan donne également de nouvelles perspectives sur la réception des travaux de ce dernier. Elle met notamment en évidence la circulation de procédés de réductions des représentations analytiques des substitutions entre différents domaines et différents auteurs. 
 
Comme nous allons l'illustrer rapidement, les travaux de Jordan en théorie des groupes ne peuvent en effet être dissociés d'autres domaines comme, par exemple, les équations différentielles. 
Nous avons déjà évoqué que Jordan avait consacré de nombreux travaux à l'étude des sous-groupes résolubles du groupe linéaire. Le mémoire "Sur la résolution algébrique des équations primitives de degré $p^2$", publié en 1868 dans le Journal de Liouville, visait ainsi à étudier les sous-groupes résolubles des groupes linéaires à deux variables ($Gl_2(F_p)$), une problématique que l'auteur présentait comme dépassant les explorations de Galois dans son "Fragment d'un second mémoire" (\hyperref[Annexe3]{\textcolor{blue}{annexe 3}}). La décomposition de la forme analytique des substitutions linéaires y jouait un rôle crucial. C'est notamment dans ce cadre que Jordan énonçait son célèbre théorème de réduction canonique des substitutions linéaires "à une forme aussi simple que possible". 

En effet, contrairement au cas des substitutions d'une variable $(x \ ax+b)$ pouvant se réduire à une combinaison d'opérations $(x \ x+1)$ et $(x \ ax)$, l'étude de substitutions de $n$ variables nécessite la prise en compte de situations plus complexes. Jordan montrait ainsi que  toute substitution linéaire à deux variables peut s'écrire sous l'une des trois formes suivantes, selon qu'une équation du second degré (l'équation caractéristique) admette deux racines réelles distinctes $\alpha, \beta$, deux racines imaginaires distinctes $\alpha + \beta i, \alpha + \beta i^p$ (où $i^2 \equiv1mod(p)$) ou une racine double :
\[
\begin{vmatrix}
z & \alpha z \\
u & \beta u 
\end{vmatrix}
, \begin{vmatrix}
z & (\alpha + \beta i)z \\
u & (\alpha + \beta i^p)u 
\end{vmatrix}
, \begin{vmatrix}
z & \alpha z \\
u & \beta z+ \gamma u 
\end{vmatrix}
\]
Ce résultat met à nouveau en évidence le rôle de modèle joué par la décomposition de la forme linéaire en deux formes de représentations analytiques de cycles : les formes canoniques réduisent les substitutions linéaires à des multiplications (deux premiers cas) ou des combinaisons multiplication-addition (troisième cas). 
Dans son traité, Jordan généralisait ce théorème à \textcolor{red}{la réduction canonique des substitutions linéaires de $n$ variables :}
\begin{quote}
\textcolor{red}{\textbf{Théorème de réduction canonique} \\ Cette forme simple
\[
\begin{vmatrix}
y_0, z_0, u_0, ..., y'_0, ... & K_0y_0, K_0(z_0+y_0), ... , K_0y'_0 \\
y_1, z_1, u_1, ..., y'_1, ... & K_1y_1, K_1(z_1+y_1), ... , K_1y'_1 \\
.... & ... \\
v_0, ... & K'_0v_0, ... \\
... & ... \\
\end{vmatrix}
\]
à laquelle on peut ramener la substitution $A$ par un choix d'indice convenable, sera pour nous sa forme canonique.\cite[p.127]{6}}
\end{quote}
Plus encore, Jordan investissait immédiatement ce résultat à d'autres domaines comme celui des équations différentielles linéaires (à coefficients constants d'abord puis dans le cas des équations de Fuchs). Dans ce cas, les substitutions présentent pourtant une nature différente, puisqu'agissant sur le corps infini des nombres complexes et non sur des corps finis. Jordan s'appuyait cependant sur la permanence de la forme analytique des substitutions afin de transférer ses procédés de réductions des groupes finis aux \textcolor{red}{systèmes d'équations différentielles}: 
\begin{quote}
\textcolor{red}{
\[
\begin{matrix}
\frac{dx_1}{dt}=a_1x_1+...+l_1x_n \\
\frac{dx_2}{dt}=a_2x_1+...+l_2x_n \\
... \\
\frac{dx_n}{dt}=a_nx_1+...+l_nx_n
\end{matrix}
\]
Ce problème peut se résoudre très simplement par un procédé identique à celui dont nous nous sommes servi, dans notre Traité des substitutions, pour ramener une substitution linéaire quelconque à sa forme canonique. Nous allons ramener de même le système à une forme canonique qui puisse s'intégrer immédiatement.
\[
\frac{dy_1}{dt}=\sigma y_1, \frac{dz_1}{dt}=\sigma z_1+y, \frac{du_1}{dt}=\sigma u_1+z_1, ..., \frac{dw_1}{dt}=\sigma w_1+v_1
\] 
[...] le système des équations aura évidemment pour intégrales le système suivant : 
\[
w_1=e^{\sigma t}\psi(t), v_1=e^{\sigma t}\psi'(t), ..., y_1= e^{\sigma t}\psi^r(t),
\] 
$\psi(t)$ étant une fonction entière arbitraire du degré $r-1$.}\cite[p.787]{7}
\end{quote}

L'étude de la circulation de tels procédés de réductions de la représentation analytique des substitutions dévoile des héritages des travaux de Jordan durant la fin du XIXe siècle. Elle permet notamment de mettre à jour l'importance de ces travaux pour l'élaboration par Henri Poincaré de la théorie des fonctions fuchsiennes en entremêlant théorie des groupes, équations différentielles, arithmétique, mécanique etc.\footnote{BRECHENMACHER (Frédéric), \href{http://hal.archives-ouvertes.fr/aut/Frederic+Brechenmacher/}{Autour de pratiques algébriques de Poincaré}, prépublication.}  Plus tard, au tournant du siècle, les procédés de réduction de Jordan allaient constituer la base d'une culture algébrique spécifique partagée par un groupe de mathématiciens français et américains. \footnote{BRECHENMACHER (Frédéric) \href{http://hal.archives-ouvertes.fr/aut/Frederic+Brechenmacher/}{Linear groups in Galois fields. A case study of tacit circulation of explicit knowledge}, Oberwolfach Reports, 2012.}

\section{Conclusion}
Les portraits de Camille Jordan ont souvent brossé les traits d'un mathématicien solitaire dans l'ombre d'Évariste Galois. Nous avons vu qu'un tel éclairage projette cependant sur le passé des organisations du savoir mathématique utilisées de nos jours : théorie des groupes, théorie de Galois, algèbre, etc. Envisagés dans la perspective de la théorie de l'ordre, les premiers travaux de Jordan nous apparaissent au contraire prendre place dès les années 1860 dans un cadre collectif. Plus encore, leurs héritages se manifestent dans des domaines variés des sciences mathématiques de la fin du XIX\up{e} siècle.

Le problème de la restitution des dimensions collectives dans lesquelles saisir la création mathématique individuelle nous invite ainsi à redécouvrir des organisations des savoirs antérieures aux disciplines qui nous contemporaines et qui, comme l'algèbre linéaire, sont centrées sur des objets comme les groupes, corps et espaces vectoriels.

Le \textit{Traité des substitutions et des équations algébriques} de 1870 présentait cependant déjà une nature hybride entre l'approche transversale de la théorie de l'ordre et une théorie centrée sur des objets. À partir de 1867, Jordan avait en effet progressivement attribué à la notion de groupe de substitutions l'"essence" qu'il avait initialement associé à la théorie de l'ordre. Il avait alors réorganisé ses travaux dans un nouveau cadre, centré non plus sur des \textit{relations transversales} mais sur des \textit{objets} comme les substitutions et les irrationnelles des équations algébriques.

Cette hybridité se manifeste notamment dans le choix de Jordan de ne pas aborder dans son \textit{Traité} les questions relatives aux polyèdres, à la cristallographie ou à la cinématique des corps solides mais, au contraire, de placer au coeur de son ouvrage une présentation des principes généraux de Galois sur les équations algébriques :
\begin{quote}
Le but de cet Ouvrage est de développer les méthodes de Galois et de les constituer en corps de doctrine, en montrant avec quelle facilité elles permettent de résoudre tous les principaux problèmes de la théorie des équations.\cite[p.VII]{6}
\end{quote}
Le Livre III, consacré aux irrationnelles, se plaçait dans un cadre collectif antérieur aux travaux de Jordan mais différent de celui de la théorie de l'ordre. Il s'inscrivait en effet dans l'héritage d'une interprétation préexistante des "principes de Galois" en termes d'"ordres d'irrationalités" définis par les racines de types d'équations algébriques non résolubles par radicaux :
\begin{quote}
De ce point de vue élevé, le problème de la résolution par radicaux, qui naguère encore semblait former l'unique objet de la théorie des équations, n'apparaît plus que comme le premier anneau d'une longue chaîne de questions relatives aux transformations des irrationnelles et à leur classification.\cite[p. VI]{6}.
\end{quote}
Cet héritage impliquait notamment des travaux d'Hermite et Kronecker qui s'étaient appuyé constamment sur des articulations entre algèbre, arithmétique et analyse.\footnote{Voir à ce sujet : GOLDSTEIN (Catherine), \href{http://smf4.emath.fr/Publications/RevueHistoireMath/17/html/}{Charles Hermite's strolls in Galois fields}, \textit{Revue d?histoire des mathématiques}, 17(2011), p. 211-270.
GOLDSTEIN (Catherine), SCHAPPACHER (Norbert), A Book in Search of a Discipline (1801-1860) in [Goldstein, Schappacher, Schwermer, \textit{op. cit.}], p. 3-66.
GOLDSTEIN (Catherine), SCHAPPACHER (Norbert), Several Disciplines and a Book (1860?1901), in [Goldstein, Schappacher, Schwermer, \textit{op. cit.}], p. 67-104.} Il s'agissait de déterminer les fonctions analytiques les plus générales permettant d'exprimer les racines d'équations de degré supérieur à cinq. Jordan mêlait cet héritage à la récente approche géométrique de Clebsch et en présentait une synthèse originale centrée sur la notion de groupe (\hyperref[Annexe4]{\textcolor{blue}{annexe 4}}).

Le \textit{Traité} de 1870 présente donc une tension entre, d'une part, sa structure sous-tendue par les procédés de décompositions analytiques des substitutions hérités du cadre transversal de la théorie de l'ordre et, d'autre part, la mise en avant d'une théorie unificatrice basée sur deux objets : les substitutions et équations algébriques. Ces deux facettes ont donné lieu à des réceptions différentes. 
La première a marqué en profondeur les pratiques algébriques de mathématiciens français puis américains. 
Certains aspects de la seconde ont dans un premier temps été envisagés dans la continuité de travaux antérieurs à ceux de Jordan, notamment ceux d'Hermite et Clebsch, tandis que l'ambition unificatrice de Jordan a été ignorée ou même vivement critiquée, notamment par Kronecker.\footnote{Voir à ce sujet : BRECHENMACHER (Frédéric), \href{http://smf4.emath.fr/Publications/RevueHistoireMath/17/html/}{La controverse de 1874 entre Camille Jordan et Leopold Kronecker}, \textit{Revue d'Histoire des Mathématiques},  13 (2007), p. 187-257.}  Plusieurs décennies plus tard, cette ambition allait au contraire être saluée comme marquant l'une des principales étapes du développement de la théorie des groupes comme domaine autonome.

\appendix
\section{Annexe 1 \\ Réceptions des travaux de Jordan \\ entre 1870 et 1914}
\label{Annexe1}

Afin de rendre compte des réceptions des travaux de Jordan, il est nécessaire de ne pas forcer ces derniers à entrer dans des cadres collectifs fixés rétrospectivement et \textit{a priori} (comme la théorie de groupes ou l'algèbre). Nous allons par conséquent rechercher des traces des lectures du \textit{Traité} chez les mathématiciens contemporains de Jordan. L'étude d'un périodique de recensions mathématiques, le \textit{\href{http://www.emis.de/MATH/JFM/JFM.html}{Jahrbuch über die Forschritte der Mathematik}}, permet de pister les publications faisant référence aux travaux de notre héros sur la période 1870-1914.\footnote{BRECHENMACHER (Frédéric), \href{ttp://hal.archives-ouvertes.fr/aut/Frederic+Brechenmacher/ }{On Jordan's measurements}, prépublication.} 

Il faut tout d'abord constater que ces publications sont classées dans des rubriques variées des classifications mathématiques de l'époque. Ainsi, alors que la rubrique "théorie des équations" s'avère très marginale et que la rubrique "théorie des substitutions" ne regroupe, avant 1895, qu'environ 10 \% des recensions faisant référence à Jordan, les rubriques "théorie des fonctions", "géométrie pure" et "géométrie analytique" sont fortement représentées. Par ailleurs, si la "théorie des groupes" devient très largement majoritaire après 1895, les références à Jordan dans ce contexte sont non seulement indépendantes de la théorie de Galois mais ce phénomène s'avère également limité dans le temps : après 1910, la grande majorité des références à Jordan portent sur le théorème de la courbe de Jordan énoncé dans le \textit{Cours d'Analyse de l'École polytechnique}.\\
Quatre principaux types de références au Traité peuvent être distingués afin de présenter un panorama simplifié des principaux usages de cet ouvrage sur la période 1870-1914. 

\subsection{Une synthèse sur la théorie des substitutions}

Ce type de références apparaît à un niveau européen peu après la parution du \textit{Traité}. Il perdure sur toute la période et s'étend aux États-Unis dans les années 1880.\footnote{Janni, Sardi, Netto, Klein, Pellet, Bolza, Hölder, Borel et Drach, Vogt, Weber, Picard, Echegaray, Bianchi, Pierpont etc.} Dans ce contexte, l'ouvrage de Jordan n'est pas envisagé comme une rupture mais au contraire dans la continuité de travaux antérieurs de mathématiciens comme Cauchy, Kronecker, Hermite, Bertrand, Serret et Mathieu. 
Les textes qui s'appuient sur des innovations spécifiques de Jordan en théorie de groupes se réfèrent quant à eux le plus souvent à des articles publiés après le \textit{Traité}. Nous ne détaillerons donc pas ces références ici.\footnote{Comme les travaux de Bochert et Maillet sur les théorèmes de finitude des groupes primitifs et transitifs ou les très nombreuses références aux travaux de Jordan de la fin des années 1870 sur la classification des groupes linéaires finis en lien avec les formes quadratiques.}  

\subsection{Les équations spéciales des fonctions elliptiques et abéliennes}
Comme le précédent, ce type de références apparaît dès 1870 à un niveau européen. Mais les publications qui l'emploient réagissent cette fois à des innovations qu'elles attribuent spécifiquement à Jordan. Les problèmes traités sont liés à la résolubilité des équations mais se présentent de manière transversale à de nombreuses rubriques comme la "théorie des formes", la "théorie des substitutions", la "géométrie analytique", les "fonctions spéciales" etc.

 Rappelons que la démonstration par Abel de l'impossibilité de résoudre par radicaux les équations algébriques générales de degré supérieur ou égal à cinq n'avait le plus souvent pas été envisagée comme la conclusion d'une longue histoire devant ouvrir l'algèbre sur les nouvelles perspectives de la théorie de Galois. De nombreux mathématiciens avaient au contraire généralisé le problème, traditionnel, de l'expression des racines des équations générales de degré inférieur à quatre par des fonctions algébriques comprenant des radicaux. Il s'agissait alors de  rechercher des fonctions analytiques les plus simples permettant d'exprimer les racines d'équations de plus haut degré. 
 
Des mathématiciens comme Hermite, Kronecker, Betti et Brioschi avaient notamment réduit le problème de la résolution de l'équation générale de degré cinq à celui de l'équation modulaire provenant du problème de la division des périodes des fonctions elliptiques, fonctions doublement périodiques sur le plan complexe.\footnote{Cette équation avait notamment été étudiée par Galois. Voir à ce sujet : GOLDSTEIN (Catherine), \href{http://smf4.emath.fr/Publications/RevueHistoireMath/17/html/}{Charles Hermite's strolls in Galois fields}, Revue d'histoire des mathématiques, 17 (2011), p. 211-270.}  Les racines de l'équation générale de degré cinq avaient ainsi été exprimées analytiquement à l'aide des fonctions elliptiques.

 Or, dans son \textit{Traité}, Jordan avait non seulement proposé une nouvelle approche des travaux d'Hermite mais il avait également énoncé l'impossibilité de résoudre les équations générales de degré supérieur à cinq par les fonctions elliptiques. Ce théorème avait été immédiatement reconnu comme l'un des résultats majeurs du \textit{Traité}.\footnote{Netto, Sylow, Marie, Kronecker, Krause, Nöther etc.}
 
 Mais Jordan avait aussi abordé la question de la recherche de fonctions permettant la résolution d'équations spéciales de degré supérieur à cinq. C'était dans ce cadre qu'avait été énoncé un autre théorème reconnu comme majeur par les contemporains : le théorème des 27 droites sur une surface cubique.\footnote{Cremona, Clebsch, Geiser, Brioschi.}\\
 Jordan avait en effet cherché à généraliser le rôle joué par les équations modulaires des fonctions elliptiques à l'équation de la division des périodes de fonctions hyperelliptiques (ou abéliennes) à quatre périodes. Il avait montré que le groupe de cette dernière équation (de degré 80) contient le groupe de l'équation des 27 droites sur une surface cubique, ouvrant ainsi la voie à des approches mêlant algèbre, analyse et géométrie dans la lignée des travaux de Clebsch et Gordan. 
C'était à l'occasion de l'énoncé de ce théorème que Jordan annonçait en 1869 la parution prochaine de son \textit{Traité} :
\begin{quote}   
Tous les géomètres connaissent le fait de l'abaissement des équations modulaires pour les transformations des degrés 5, 7 et 11, et les importantes conséquences qu'en a déduites M. Hermite. MM. Clebsch et Gordan ont signalé un abaissement analogue pour les équations des périodes dont dépend la bissection des fonctions abéliennes. Nous venons d'obtenir un résultat du même genre pour l'équation qui donne la trisection dans les fonctions à quatre périodes. [...] La réduite $Z$ présente cette particularité remarquable d'avoir le même groupe que l'équation $X$ qui détermine les vingt-sept droites situées sur une surface du troisième ordre. [...] Mais cette proposition exigeant quelques développements, nous en réservons la démonstration pour le Traité des équations algébriques que nous nous occupons de publier.
\end{quote}
Le théorème des 27 droites jouait en effet un rôle clé dans la troisième section du \textit{Traité}, "Des irrationnelles" (voir à ce sujet \hyperref[Annexe4]{\textcolor{blue}{l'annexe 4}}).

\subsection{Les solides réguliers}
Ce type de références témoigne de lectures du \textit{Traité} en relation avec d'autres travaux de Jordan : polyèdres, cristallographie, cinématique, surfaces (y compris les surfaces de Riemann) ou équations différentielles.   Bien que Jordan ait choisi de ne pas inclure de tels sujets dans son \textit{Traité} de 1870, ces problèmes étaient fortement interconnectés à cette époque. 

Ces sujets mettent notamment tous en jeu des groupes de substitutions. En effet la donnée d'un polyèdre régulier définit un groupe, le groupe de symétries du solide, pouvant être envisagé d'un point de vue cinématique (mouvements d'un solide autour d'un axe par exemple) et dont l'étude peut non seulement être appliquée à la cristallographie mais aussi à l'étude des symétries des surfaces (y compris des surfaces de Riemann) et, par là, des équations différentielles. Cette approche allait notamment être développée dans les années 1870 par les travaux de Felix Klein.

Comme pour le précédent, ce type de références s'avère transversal à de nombreuses rubriques de la classification mathématique : "surfaces", "cinématique", "physique moléculaire", "fonctions spéciales" etc. 

\subsection{Groupes linéaires et champs de Galois}
Ce type de références se manifeste plus tardivement que les deux précédents et de manière plus locale : il implique essentiellement un groupe de mathématiciens français et américains sur la période 1893-1907. 
 \begin{center}
\includegraphics[scale=0.6]{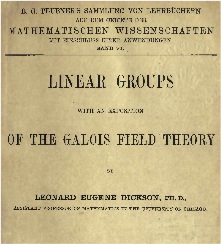}\\*
L'un des auteurs les plus importants de ce groupe, Leonard Dickson de Chicago, publiait en 1901 l'ouvrage : \textit{Linear groups with an exposition of the Galois field theory}, présentant des généralisations systématiques de résultats du \textit{Traité}.
\end{center}
Dans ce contexte, un type spécifique de lecture du \textit{Traité} est à la base d'une véritable culture algébrique commune à des mathématiciens français et américains.\footnote{BRECHENMACHER (Frédéric) \href{http://hal.archives-ouvertes.fr/aut/Frederic+Brechenmacher/}{Linear groups in Galois fields. A case study of tacit circulation of explicit knowledge}, Oberwolfach Reports, 2012.} Comme le manifeste l'expression "Jordan's linear groups" utilisée par ces derniers, Jordan était cité principalement pour un théorème qui faisait déjà l'objet de sa première thèse en 1860 et qui jouait un rôle structurant pour le traité de 1870 : l'origine des groupes linéaires finis. 

C'est aussi principalement dans ce contexte que les travaux de Jordan étaient cités en compagnie de ceux de Galois. Mais la référence à ce dernier ne portait pas sur ce que nous désignons aujourd'hui comme "théorie de Galois" mais sur ce qui était alors dénommé les "imaginaires de Galois", "champs de Galois" ou "Galois fields", c'est-à-dire les corps finis introduits par Galois dans sa note "Sur la théorie des nombres" publiée en 1830. 

Comme nous le montrons dans le corps du texte et en\hyperref[Annexe4]{ \textcolor{blue}{annexe 4}}, les imaginaires de Galois permettent d'indexer les lettres permutées par des groupes finis et ainsi de représenter les substitutions par une forme analytique. Dans les travaux de Jordan, cette forme de notation s'accompagnait de procédures spécifiques de réductions de la représentation analytique des substitutions linéaires. Dans les années 1870-1880, ces procédures avaient circulé en profondeur dans des travaux de mathématiciens comme Henri Poincaré. Ces derniers ne faisaient pas pour autant toujours explicitement référence aux travaux de Jordan. Pour cette raison, ce type de référence n'avait émergé au grand jour que dans les années 1890. Il s'agit là de l'un des principaux héritages des travaux de jeunesse de Jordan.

\subsection{Représentations simplifiées des répartitions des types de références à Jordan}

\begin{center}
\includegraphics[scale=0.6]{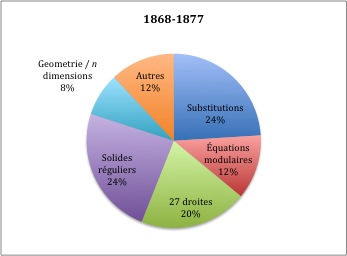}\\*
\end{center}

 \begin{center}
\includegraphics[scale=0.6]{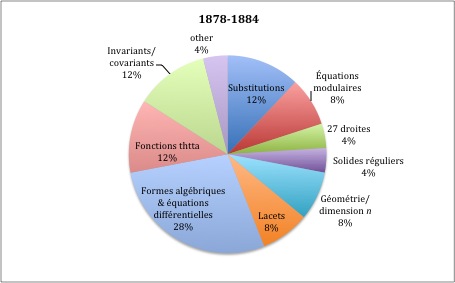}\\*
\end{center}

 \begin{center}
\includegraphics[scale=0.6]{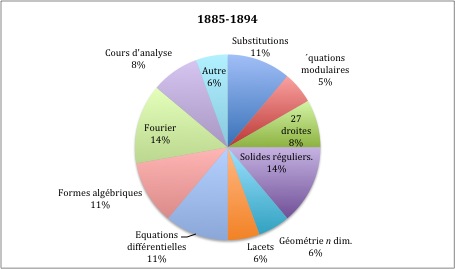}\\*
\end{center}

 \begin{center}
\includegraphics[scale=0.6]{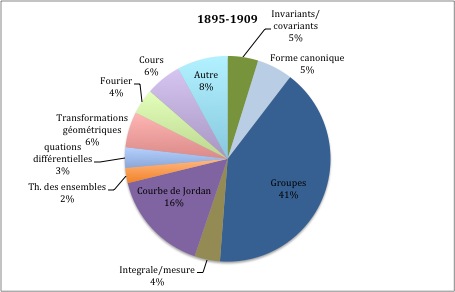}\\*
\end{center}

 \begin{center}
\includegraphics[scale=0.6]{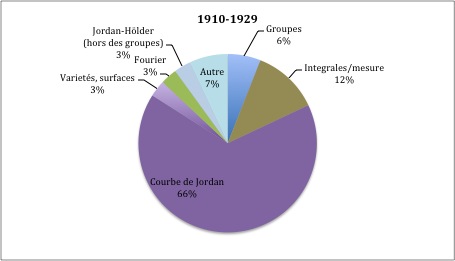}\\*
\end{center}

\section{Annexe 2. \\ Nombre des valeurs d'une fonction \\ et résolution algébrique des équations}
\label{Annexe2}
\subsection{L'équation quadratique}

Prenons pour premier exemple le cas de l'équation quadratique sur le corps des nombres rationnels :
\[
x^2-c_1x+c_2=0
\]
Les coefficients $c_1$ et $c_2$ sont des fonctions symétriques des racines $x_1$ et $x_2$. Il s'agit donc de fonctions qui ne prennent qu'une seule valeur par permutations des racines :
\[
c_1= x_1+x_2 \ ; \ c_2 = x_1x_2
\]
Toute fonction qui ne prend qu'une seule valeur peut être exprimée rationnellement dans le corps de nombres auxquels appartiennent les coefficients $c_1$ et $c_2$, c'est à dire ici sur $\Q$. Au contraire, la fonction 
\[
x_1-x_2
\]
prend deux valeurs par permutations des racines et n'est donc pas rationnellement connue. On peut alors s'intéresser au groupe des permutations qui laissent fixe cette fonction. Ce groupe laisse également  fixe la racine $x_1$ et, par conséquent, $x_1$ et $x_1-x_2$ peuvent être exprimée rationnellement l'une par rapport à l'autre : 
\[
x_1=\frac{c_1+x_1-x_2}{2}
\]
Ainsi, si l'on adjoint au corps de nombres initial le nombre $x_1-x_2$, on obtient également le nombre $x_1$ : au XIX\up{e} siècle, les fonctions de plusieurs variables, ou résolventes, étaient employées pour penser ce que nous envisageons aujourd'hui en termes de structures algébriques de corps et d'extensions de corps.\\
Revenons à présent à l'équation quadratique. Le discriminant $\Delta$ ci-dessous est une fonction d'une seule valeur et peut donc être exprimé à l'aide des coefficients $c_1$ et $c_2$ :
\[
\Delta=(x_1-x_2)^2=(x_1+x_2)^2-4x_1x_2=c_1^2-4c_2
\]
La fonction $\sqrt{\Delta}$ est, quant à elle, une fonction de deux valeurs, tout comme les racines $x_1$ et $x_2$ elles mêmes. Ces fonctions peuvent donc être exprimées rationnellement l'une par rapport à l'autre avec les formules bien connues de résolution par radicaux de l'équation quadratique.

\subsection{L'équation de troisième degré}

Dans le cas de l'équation du troisième degré 
\[
x^3-c_1x^2+c_2x-c_3=0
\]
la résolution nécessite la détermination de trois fonctions $x_1$, $x_2$, $x_3$ prenant trois valeurs  par permutations, c'est à dire finalement d'une fonction prenant $3!=6$ valeurs. Soit $a_1, a_2, a_3$ trois paramètres, la fonction $\xi = a_1x_1+a_2x_2+a_3x_3$ est une telle fonction de $3!=6$ valeurs. Il s'agit là de ce que l'on a appelé à la suite des travaux d'Enrico Betti une "résolvente de Galois". Si l'on parvient à exprimer les racines à l'aide de radicaux, alors $\xi$ est également exprimable par radicaux et réciproquement.

En un sens, le problème de la résolution algébrique de l'équation consiste à passer des trois fonctions d'une valeur $c_1, c_2, c_3$ à la fonction de six valeurs $\xi$.\\ 
Comme dans le cas de l'équation quadratique, la racine du déterminant donne une fonction de deux valeurs à partir de laquelle toute fonction de deux valeurs peut être exprimée rationnellement. Pour résoudre le problème, il suffit donc de trouver une fonction des racines dont une certaine puissance prend deux valeurs. Tel est le cas du cube de ce que l'on appelle la résolvente de Lagrange :
\[
\phi = x_1+ \omega^2 x_2+\omega x_3
\]
où $\omega$ est une racine primitive de l'unité ($\omega^3=1$). \\
On peut donc exprimer $\phi$ rationnellement à l'aide de $\sqrt \Delta$:
\[
\phi^3= \frac{1}{2}(2c_1^3-9c_1c_2+27c_2+3\sqrt(-3\Delta))
\]
et l'on déduit de l'expression ci-dessus les célèbres formules de Cardan.

\subsection{En général}

Résoudre une équation générale de degré $n$ implique la considération d'une fonction résolvente prenant $n!$ valeurs distinctes. C'est à partir de la considération de l'ensemble des substitutions laissant invariantes de telles fonctions et leurs facteurs polynomiaux que Galois définissait ce que l'on appelle aujourd'hui le groupe de Galois d'une équation.\footnote{C'est à dire, en termes contemporains, laissant stable un corps de racines.}
Nous avons vu qu'une résolvente de Galois est une fonction de $n$ variables prenant $n!$ valeurs distinctes. Elle n'est par conséquent invariante que par le groupe trivial réduit à la permutation identité. A l'opposé, une fonction symétrique est invariante par toutes les substitutions du groupe symétrique. 
Adjoindre des racines à une équation comme nous l'avons fait plus haut avec $\sqrt{\Delta}$ implique une décomposition de la résolvente de Galois en facteurs. A chacun de ces facteurs peut être associé le groupe des substitutions le laissant invariant.

Par exemple, pour $n=4$, la fonction  suivante peut être interprétée comme exprimant les relations entre les racines d'une équation irréductible du quatrième degré :
\[
\phi= x_1x_2+x_3x_4
\]
Cette fonction prend trois valeurs pour toutes les $4!=24$ permutations de $\Sigma(4)$ :
\[
x_1x_2+x_3x_4, x_1x_3+x_2x_4, x_1x_4+x_2x_3
\]
Le groupe $G$ associé à la fonction $\phi$ est composé des huit  substitutions laissant $\phi$ globalement invariante :
\[
G ={I, (x_1x_2), (x_3x_4), (x_1x_2)(x_3x_4), (x_1x_3)(x_2x_4), (x_1x_4)(x_2x_3), (x_1x_3 x_2x_4), (x_1x_4, x_2x_3 )}.
\]

\section{Annexe 3. \\ La représentation analytique des substitutions}
\label{Annexe3}

\subsection{Représentation et indexation}

Au cours du XIX\up{e} siècle, différentes formes de représentations des substitutions ont été développées. Citons notamment:
\begin{itemize}
\item la représentation en deux lignes ($a$ s'envoie sur $d$ ; $b$ sur $c$ etc.) :
\[
(a, b, c, d, e, ...)
(d, c, a, e, b, ...)
\]
\item la représentation par produits de transpositions :
\[
(a d) (d e) (eb) (bc) (c a)
\]
\item la représentation symbolique des opérations sur les substitutions :
\[
ghg^{-1}=k
\] 
\item la représentation tabulaire, notamment utilisée par Galois comme l'a analysé Caroline Ehrhardt\footnote{EHRHARDT (Caroline), \textit{Evariste Galois. La fabrication d'une icône des mathématiques}, Paris, Éditions de l'EHESS, 2011.}
\begin{center}
\includegraphics[scale=0.7]{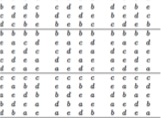}\\*
\end{center}
\item la représentation analytique. 
\end{itemize}
Cette dernière forme de représentation a joué un rôle important bien qu'elle soit passée inaperçue de nombreux travaux historiques. Elle consiste à indexer les lettres par des entiers afin de représenter les substitutions par des polynômes. 

Étant donné une  substitution $S$ sur $p$ lettres $a_0, a_1, ..., a_{p-1}$, le problème consiste à déterminer une fonction analytique $\phi$ telle que $S(a_i)=a_{\phi(i)}$. La résolution peut s'appuyer sur l'utilisation de polynômes d'interpolation de Lagrange. Par exemple, pour $p=5$, Hermite a démontré que les formes réduites de telles substitutions sont de la forme :
\[
i \ ; i^2 \ ; i^3+ai.
\]

Dans le cas d'un nombre premier $p$ de lettres, une indexation peut être obtenue par cyclotomie en considérant ces lettres comme des racines $p^e$ de l'unité, c'est à dire comme des racines de l'équation binôme $x^p=1$ et dont toutes les $p$ racines sont obtenues par la suite des puissances d'une racine primitive $\omega$  de l'unité : 
\[
\omega^0, \omega^1, ... , \omega^i, ..., \omega^{p-1}
\]

Cette indexation, $i=0, 1, 2, ..., p-1$, correspond donc à un groupe cyclique engendré par la substitution représentée analytiquement par $(i \ i+1)$.
Mais les racines peuvent aussi être réindexées en considérant une racine primitive $g$ de l'équation cyclotomique $mod.p$
 \[
x^p \equiv 1 \ mod(p)
\]

afin d'obtenir la suite :
\[
\omega^g, \omega^{g^2}, ..., \omega^{g^{p-1}}
\]
correspondant au même groupe cyclique mais dont la substitution génératrice est cette fois représentée sous la forme analytique $(i \ gi)$. \\
Dans ses célèbres \textit{Disquitiones arithmeticae} de 1801, Gauss avait étudié la résolubilité par radicaux des équations cyclotomiques, c'est à dire des équations irréductibles déduites des équations binômes :\footnote{NEUMANN (Olaf), The Disquisitiones Arithmeticae and the Theory of Equations, in GOLDSTEIN (Catherine), SCHAPPACHER (Norbert), SCHWERMER (Joaquim) (eds.),  \textit{The Shaping of Arithmetics after C. F. Gauss's Disquisitiones Arithmeticae}, Berlin: Springer, 2007, p. 107-128.}  
 \[
\frac{x^p-1}{x-1}=x^{p-1}+x^{p-2}+ ...+x+1
\]
A cette occasion, l'ensemble des racines d'une telle équation avait été décomposé en ce que Jordan désignait comme des systèmes imprimitifs. Cette décomposition s'appuyait sur des réindexations des racines en un ordre spécifique à l'aide des indexations données par une racine primitive de l'unité $\omega$ et d'une racine primitive $g$ de l'unité modulo $p$, c'est à dire par les deux formes de représentations analytiques des cycles. \\
Détaillons un exemple pour $p=19$. Dans ce cas, $p-1=18$ peut se factoriser en $18 =3.6$. Par conséquent, les 18 racines cyclotomiques
\[
\omega, \omega^{2}, ..., \omega^{18}
\]
peuvent se répartir en 6 "périodes" (pour reprendre le terme de Gauss) de 3 racines et donc en un produit d'équations de degrés 3 et 6. Par exemple, la période des 3 racines $\eta_1, \eta_2, \eta_3$ de l'équation 
\[
x^3+x^2-6x-7=0
\]
correspond aux sommes suivantes :
\[
\eta_1=\omega+\omega^7+\omega^8+\omega^11 +\omega^12 +\omega^18 
\eta_2= \omega^2+\omega^3+\omega^5+\omega^14+\omega^16+\omega^17
\eta_3= \omega^4+\omega^6+\omega^9+\omega^10+\omega^13+\omega^15
\]
Les exposants de chaque suite de termes impliqués dans ces sommes  sont indexés par des puissances successives d?une racine primitive $mod. 19$. Ici, par exemple, $g=2$ est une racine primitive de l'unité $mod. 19$ (car $218 = 262 144 = 1 + 19. 13 797$) et les indexations des  $\eta_i$ correspondent aux 3 cycles des 6 puissances de $2^3 \ mod. 19$ : 
\begin{align*}
(2^3)^0 \equiv1 mod (19) \\
(2^3)^1 \equiv 8 mod (19) \\
(2^3)^2 \equiv 7 mod (19) \\
(2^3)^3 \equiv 18 mod (19) \\
(2^3)^4 \equiv 11 mod (19) \\
(2^3)^5 \equiv 12 mod (19)
\end{align*}
Soit, finalement :
\[
1, 7, 8, 11, 12, 18
\]

Les exposants des deux autres suites sont alors obtenus en multipliant cette suite par $ g \equiv 2  \ mod (19)$ puis $g \equiv 2^2 \ mod (19)$. \\
En général, pour toute factorisation de $p-1=ef$, posons $h=ge$ et considérons l'équation de degré $e$ dont les racines correspondent aux $e$ périodes de sommes de $f$ termes: 
\[
\eta_i= \omega^i+omega^{ih}+ ...+\omega^{ih^{f-1}} (1\leq i \leq e)
\]

Cette décomposition des racines en périodes donne une factorisation de l'équation initiale en $e$ facteurs de degrés $f$.

Comme l'a analysé récemment Jenny Boucard, Poinsot avait insisté, lors d'un commentaire des travaux de Gauss en 1808, sur le procédé de décomposition simultané en "groupes"  des substitutions et des lettres sur lesquelles agissent ces substitutions.\footnote{Voir à ce sujet : BOUCARD (Jenny), Louis Poinsot et la théorie de l'ordre : un chaînon manquant entre Gauss et Galois ?, \textit{Revue d?histoire des mathématiques}, 17, fasc. 1 (2011), p. 41-138.
FREI (Gunther), The unpublished section eight: On the way to function fields over a finite field, in GOLDSTEIN (Catherine), SCHAPPACHER (Norbert), SCHWERMER (Joaquim) (eds.), \textit{The Shaping of Arithmetics after C. F. Gauss's Disquisitiones Arithmeticae}, Berlin: Springer, 2007, p. 159-198.
NEUMANN (Olaf), The Disquisitiones Arithmeticae and the Theory of Equations, in GOLDSTEIN (Catherine), SCHAPPACHER (Norbert), SCHWERMER (Joaquim) (eds.), \textit{op. cit.}, p. 107-128.}
  Il avait aussi à cette occasion mis en avant les deux formes de représentation analytique des cycles : dans la méthode de décomposition de Gauss, les "groupes" de racines peuvent être représentées "comme sur un cercle" de sorte que l'on passe d'une racine à l'autre en "avançant" d'un rang dans la liste des index $i$ des puissances de la racine primitive $\omega$,  c'est à dire par l'opération d'addition $(i \  i+1)$ tandis que l'on passe d'un "groupe" à l'autre par "rotation" du cercle sur lui même, c'est à dire par l'opération de multiplication $(i \  gi)$. 
\begin{center}
\includegraphics[scale=0.7]{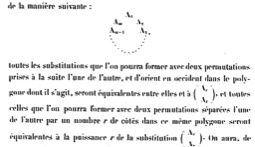}\\*
\end{center}

Un même type de représentation circulaire avait été utilisé par Cauchy dans ses travaux sur les nombres de valeurs des fonctions publiés en 1815. Contrairement à Poinsot, Cauchy avait cependant préféré d'autres modes de représentations à la représentation analytique (notation symbolique, produits de transpositions, représentations tabulaires).\footnote{Voir DAHAN (Amy), Les travaux de Cauchy sur les substitutions. Étude de son approche du concept de groupe, \textit{Archive for History of Exact Sciences}, vol. 23, 1980, p. 279-319.} Il avait introduit les termes de "substitutions arithmétiques" pour le cas $(i \ i+a)$ et "substitutions géométriques" pour le cas $(i \ gi)$. Ces termes avaient été repris par la grande majorité de présentations de la théorie des substitutions jusqu'à la fin du XIX\up{e} siècle. Les travaux de Jordan revendiquaient quant à eux l'héritage de Poinsot et développaient des considérations générales sur les groupes imprimitifs sur $p^n$ lettres. Décomposés sur le modèle de la méthode de Gauss, ces groupes engendrait des "groupes linéaires" de substitutions de forme analytique $(i \ ai+b)$.

\subsection{Le critère de résolubilité de Galois}

La représentation analytique des substitutions ne se limitait en effet pas à une simple notation mais s'accompagnait de procédures spécifiques de décompositions sur le modèle de la méthode de Gauss. Cette situation est bien illustrée par le rôle clé joué par cette représentation dans les travaux de Galois.

Rappelons que le théorème qui concluait le célèbre "Mémoire sur les conditions de résolubilité des équations par radicaux" énonçait un critère de résolubilité d'une classe d'équations algébriques : les équations dont toutes la racines sont fonctions rationnelles de deux d'entre elles. 

Cette classe d'équations généralisait les équations étudiées par Abel, pour lesquelles toutes les racines sont fonctions d'une seule d'entre elles, équations généralisant elles mêmes les équations cyclotomiques de Gauss pour lesquelles toutes les racines peuvent être obtenues comme puissances successives d'une racine primitive.

De fait, la liste de propositions donnée par Galois ne s'arrêtait pas à l'exposé général de la première partie du Mémoire dans laquelle se trouve la célèbre présentation du problème de la résolubilité en termes d'adjonctions de racines et décompositions successives de groupes. Au contraire, la partie d'applications du mémoire visait notamment à donner un critère de résolubilité pour les équations irréductibles de degré premier.

La représentation analytique des substitutions jouait un rôle crucial dans l'énoncé et la démonstration du théorème de Galois. Ainsi, pour qu'une équation irréductible de degré premier soit résoluble par radicaux :
\begin{quote}
 \textit{il faut} et \textit{il suffit} que toute fonction invariable par les substitutions
 \[
x_i, x_{ai+b}
\]
soit rationnellement connue.\cite[p.431]{11}
\end{quote}

Dans le cadre de la théorie de Galois qui nous est contemporaine, ce théorème peut s'interpréter comme énonçant qu'une équation irréductible de degré $p$ est résoluble par radicaux si et seulement si son groupe de Galois est un sous groupe du groupe affine.  Mais cette interprétation rétrospective ne permet pas de percevoir le rôle crucial joué par la représentation analytique des substitutions. Ce rôle se manifeste au contraire dans la démonstration de Galois qui consistait à rechercher la forme analytique f des substitutions $(x_i, x_{f(i)})$ les plus générales \textcolor{blue}{transformant un cycle $(i \ i+a)$ en un cycle},\footnote{\textcolor{blue}{Il s'agit en termes contemporains de rechercher le plus grand groupe dans lequel un $p$ groupe abélien élémentaire est distingué. Comme nous le commentons dans le coeur du texte, il s'agit également de la manière dont Jordan introduit le groupe linéaire dans son \textit{Traité} de 1870.}}  c'est à dire telle que
 \[
f(i+a)=f(i)+A
\]

 Par conséquent
 \[
f(i+2a)=f(i)+2A,..., f(i+ma)=f(i)+mA
\]
 Si $a=1$ et $i=0$, alors $f(i)=ai+b$.
Galois désignait de telles substitutions par le terme " substitution linéaire".\\
Le critère de résolubilité de Galois allait constituer la principale référence aux travaux de ce dernier jusqu'à la fin du XIX\up{e} siècle. C'était notamment ce critère qui était mis en avant par Liouville lors de l'édition des oeuvres de Galois en 1846. C'était aussi par ce critère que se concluaient toutes les présentations des travaux de Galois dans les traités publiés jusqu'à la fin du siècle, à l'exception de ceux de Jordan et Klein. 

\subsection{Le Fragment de second mémoire}

Contrairement à la majorité de ses contemporains, Jordan s'était davantage intéressé au "Fragment d'un second mémoire" de Galois qu'au célèbre "Mémoire sur la résolubilité des équations algébriques".\footnote{Rappelons que si ce texte n'a été publié qu'en 1846, il n'en est probablement pas moins antérieur à la version finale du premier mémoire en raison des épisodes bien connu des pertes successives des deux premières versions soumises par Galois à l'Académie.} Dans ce fragment, Galois avait tenté de généraliser son théorème aux équations irréductibles de degré $p^n$. Afin d'indexer des systèmes de $p^n$ lettres, il avait partagé les lettres du "groupe", c'est à dire ici de l'ensemble des racines de l'équation,\footnote{Rappelons que chez Galois, le terme "groupe"  pouvait selon les cas désigner des systèmes de substitutions et des blocs de lettres. Cette ambiguïté a souvent été soulignée par les historiens. Remarquons cependant qu'il s'agit là de la nature même du concept de groupe tel qu'il se présente par exemple dans les périodes de Gauss.} en $n$ groupes d'un même nombre de lettres $p$, représentés sous la forme tabulaire suivante :
 \[
\begin{matrix}
a_0 & a_1 & a_2 & ... & a_{p-1} \\
b_0 & b_1 & b_2 & ... & b_{p-1} \\
c_0 & c_1 & c_2 & ... & c_{p-1}  \\
.. & .. & .. & .. & ..
\end{matrix}
\]
Cette décomposition  permettait une indexation des lettres en $n$ séries de $p$ indices :
\begin{quote}
La forme générale des lettres sera
 \[
a_{\begin{matrix}
k,&k,&k,&...&k,\\
1&2&3&...&\mu
\end{matrix}}
\]
$\begin{matrix}
k, & k, & k, & .... & k,\\
1 & 2 & 3 & ... & \mu
\end{matrix}$
étant des indices qui peuvent prendre chacun les $P$ valeurs $0, 1, 2, 3, ..., P-1$.
\end{quote}
L'indexation visait à introduire une représentation analytique des substitutions de $p^n$ lettres par des fonctions $\phi, \psi, \chi, ... \sigma$ des \textcolor{blue}{indices} :\footnote{\textcolor{blue}{En termes contemporains, les indices forment un corps fini de $p^n$ éléments défini comme un espace vectoriel de dimension $n$ sur le corps $F_p$}}.  

\begin{quote}
[...]  dans le groupe $H$, toutes les substitutions seront de la forme
 \[
\begin{bmatrix}
a_{\begin{matrix}
k,&k,&k,&...&k,\\
1&2&3&...&\mu
\end{matrix}} & 
a_{\begin{matrix}
\phi(k),&\psi(k),&\chi(k),&...&\sigma(k),\\
1&2&3&...&\mu
\end{matrix}}
\end{bmatrix}
\]
\end{quote}
Comme pour le cas des équations de degré premier, Galois recherchait alors la représentation analytique des substitutions \textcolor{blue}{"transformant les cycles en cycles"}.\footnote{\textcolor{blue}{En termes contemporains, il s'agit de rechercher le groupe maximal dans lequel un groupe abélien élémentaire est un sous groupe distingué.}} Il ne donnait cependant une solution de ce problème que pour le cas des équations primitives de degré $p^2$ .  Dans ce cas, les cycles prennent la forme 
 \[
\begin{bmatrix}
a_{\begin{matrix}
k,&k',\\
1&2
\end{matrix}} & 
a_{\begin{matrix}
k+\alpha,& k'+\alpha'\\
1 1& 2 2
\end{matrix}}
\end{bmatrix}
\]
et Galois montrait que la forme analytique des substitutions cherchées est une "forme linéaire" : 
 \[
\begin{bmatrix}
a_{\begin{matrix}
k,&k',\\
1&2
\end{matrix}} & 
a_{\begin{matrix}
mk+n,& mk'+n\\
1 1& 2 2
\end{matrix}}
\end{bmatrix}
\]
Comme nous le détaillons dans le corps de l'article, Jordan avait présenté d'une manière semblable  l'"origine du groupe linéaire" sur $p^n$ lettres dans son \textit{Traité} de 1870.

\subsection{Les imaginaires de Galois}
La représentation analytique des substitutions était également liée à l'introduction par Galois des "imaginaires de la théorie des nombres" dans une note, publiée dans le \textit{Bulletin de Férussac} en juillet 1830, et présentée comme un lemme pour l'étude des substitutions primitives. \footnote{Avant Galois, les équations de congruences et corps finis avaient été considéré par Gauss  ainsi que par Poinsot. Voir Frei, 2007, \textit{op. cit} et Boucard, 2011, \textit{op. cit}.} \\
Nous avons vu que dans le cas de substitutions permutant un nombre premier de lettres, la suite d'indices $0, 1, 2, ..., p-1$ peut être réindexée par la suite $1, g , g^2, ...,g^{p-1}$ ou $g$ est une racine primitive de l'équation $x^p \equiv 1 mod (p)$. Nous avons vu également que de telles réindexations permettent de passer d'une forme d'action des cycles à une autre, c'est à dire de la représentation analytique $(i \ i+1)$  à $(i \ gi)$. 
Dans sa "Note sur la théorie des nombres", Galois avait généralisé ce mode d'indexation au cas de $p^n$ lettres. Il considérait les $p^n$ expressions algébriques  
 \[
a^{j^{n-1}}+b^{j^{n-2}}+...+1\\
\]
formées avec une racine imaginaire $j$ d'une équation irréductible de degré $n$ modulo $p$, 
 \[
f(x) \equiv 0 \ mod (p)
\]
et montrait que des telles expressions peuvent être considérées comme les $p^n$ racines de l'équation
 \[
x^{p^n} \equiv x \ mod (p)
\]
Un système de $p^n$ indices peut être ainsi \textcolor{blue}{réindexé}\footnote{\textcolor{blue}{En termes contemporains, un corps fini de $p^n$ éléments définit à la fois un groupe additif, ou un espace vectoriel de dimension $n$ sur $F_p$, et un groupe cyclique de $p^{n-1}$ éléments, son groupe multiplicatif.}} par analogie avec l'indexation cyclotomique de $p$ lettres par la suite des puissances d'une racine primitive $j$ de la congruence ci-dessus. \\
La conclusion de la note est consacrée à l'utilisation d'une telle indexation pour les racines d'une équation primitive de degré $p^n$. Une représentation analytique des substitutions sur $p^n$ lettres est ainsi obtenue par le recours à un seul paramètre $j$, à la différence des suites d'indices employées dans le fragment de second mémoire.

Galois en déduisait un critère de résolubilité des équations primitives de degré $p^n$ que Jordan allait plus tard contredire, à savoir que de telles équations sont résolubles si et seulement si leurs substitutions prennent une forme linéaire. Il n'avait cependant pas donné de démonstration de son énoncé et s'était contenté d'indiquer que celui-ci découlait de la décomposition de la forme linéaire des substitutions 
 \[
(j (aj+b)^{p^r })
\]
en un produit
 \[
a'(i+b')^{p^r }
\]
sur le modèle de la démonstration du critère de résolubilité des équations de degré $p$ : "les personnes habituées à la théorie des équations le verront sans peine".\cite[p.435]{11}

\subsection{Principes généraux et applications}

Il est bien connu que Galois avait distingué entre les principes généraux qu'il avait exposé dans son Mémoire et leurs applications aux équations de degré premier, aux équations primitives de degré $p^n$ et aux équations modulaires.

Durant le XX\up{e} siècle, les commentaires sur Galois se sont souvent focalisés sur ces principes généraux dans le cadre de quêtes d'origines de la théorie de Galois et de la théorie des groupes. Il faut cependant être attentif au rôle de modèle joué par les "applications". La décomposition de la forme analytique des substitutions linéaires donne en particulier un modèle aux principes plus généraux de décompositions des groupes. Ce rôle de modèle se manifeste aussi dans le \textit{Traité} de Jordan en 1870 comme nous le décrivons dans le corps du texte.

Dans sa célèbre lettre à Chevalier, Galois lui même soulignait le rôle de modèle joué par la "méthode de décomposition de M. Gauss" pour un premier type de décomposition des groupes. Déjà dans sa première note sur la question de la résolubilité des équations - parue en avril 1830 dans le \textit{Bulletin de Ferrusac} - Galois avait introduit une distinction entre \textcolor{blue}{équations primitives et non primitives}\footnote{\textcolor{blue}{Une équation non primitive de degré $mn$ est une équation qui peut être décomposée en m facteurs de degrés $n$ par le moyen d'une seule équation de $m$}} appelées également "équations de M. Gauss". Dans le "Fragment de second mémoire", le problème des équations résolubles de degré composé était précisément réduit par la "méthode de décomposition de M. Gauss" à celui des équations résolubles primitives de degré $p^n$.

Au contraire, l'étude des cas de réduction de degré des équations modulaires, associées à la forme analytique homographique $\frac{ai+b}{ci+d}$ ($ad-bc \ne 0$), avait servi de modèle à un autre mode de décomposition des groupes et avait par là mis en évidence la spécificité du mode précédent que Galois avait alors dénommé \textcolor{blue}{"décomposition propre"}.\footnote{\textcolor{blue}{Le premier cas correspond en termes contemporains à une décomposition en sous groupes distingués tandis que le second s'appuie sur des sous groupes non distingués.}}

Nous avons vu que cette décomposition propre était indissociable de celle de la forme analytique linéaire $(i \ ai+b)$ en deux formes d'actions des cycles $(i \ i+1)$ et $(i \ gi)$. C'était encore sur ce modèle que Galois avait étendu son critère de résolubilité aux équations résolubles primitives de degre $p^n$ en énonçant que ces dernières devaient être associées à la forme linéaire :  
 \[
x_{k, l, m ...} | x_{ak+bl+cm+...+h, a'k+b'k+c'm+...+h', a''k+ ...}
\]
La réfutation de cet énoncé allait constituer le point de départ de la relation établie par Jordan à Galois en 1860.

\subsection{La représentation analytique des substitutions et les réceptions des travaux de Galois}

Les trois formes de représentations analytiques des substitutions utilisées par Galois permettent également de suivre différentes formes de réceptions des travaux de ce dernier. En effet, à l'exception des commentaires systématiques d'Enrico Betti au début des années 1850 et de ceux de Camille Jordan à la fin des années 1860, les références à Galois ont le plus souvent visé une partie précise des travaux de ce dernier.
\begin{itemize}

\item La forme $(i \ ai+b)$

Comme nous l'avons vu, les substitutions linéaires à une variable étaient associées au théorème concluant le Mémoire de Galois et donnant un critère de résolubilité des équations irréductibles de degré premier. Ce théorème avait constitué la principale forme de référence aux accomplissements de Galois jusqu'à la fin du XIX\up{e} siècle, c'est à dire pendant la période séparant la mise en lumière de ce théorème par Liouville lors de la première édition des oeuvres de Galois et sa disparition complète de la nouvelle préface rédigée par Picard en 1897. \\
Le théorème de résolubilité de que l'on a appelé plus tard les "équations de Galois" concluait ainsi la grande majorité des présentations de la théorie de Galois des équations algébriques à l'exception de celles de Jordan et de Klein qui n'y faisaient pas même mention.\footnote{Ainsi en est-il des traités de Serret (1866), Netto (1882), Bolza (1891), Borel et Drach (1895), Vogt (1895), Weber (1895), Picard (1897), Pierpont (1900).}

\item La forme $(i \ \frac{ai+b}{ci+d})$

Contrairement aux traités, les articles faisant références aux travaux de Galois mentionnent rarement le critère de résolubilité et bien plus souvent l'énoncé de Galois relatifs à la réduction du degré des équations modulaires d'ordre 5, 7 et 11. L'intérêt d'Hermite, Betti et Kronecker pour cet énoncé se manifeste dès le début des années 1850 et s'affirme à la fin de la décennie avec l'utilisation de l'équation modulaire d'ordre 5 pour donner une expression analytique des racines de l'équation générale du cinquième degré à l'aide des fonctions elliptiques. \\
Dans la seconde moitié du XIX\up{e} siècle, la référence aux travaux de "Galois-Betti-Hermite" sur l'équation modulaire devient l'une des principales formes de références aux travaux de Galois. Elle intervient notamment de manière cruciale dans le \textit{Traité} de Jordan en 1870 (\hyperref[Annexe4]{\textcolor{blue}{annexe 4}}).

\item La forme linéaire générale 
 \[
 (k, l, m ... ; \  ak+bl+cm+...+h, a'k+b'k+c'm+...+h', a''k+ ...)
 \]
Les références aux travaux de Galois sur les équations de degré composé sont nettement plus rares que les précédentes et s'accompagnent d'usages de la représentation analytique des substitutions de $p^n$ variables ainsi que des imaginaires de Galois. \\
Après Betti en 1851, on trouve une telle référence dans deux courtes notes publiées en 1856 par Alexandre Allégret. Ces notes visaient à généraliser le critère de Galois aux équations de degré composé, problématique que l'auteur avait relié aux travaux récents de Kronecker, Betti et Pierre-Laurent Wantzel. Allégret avait par conséquent été l'un des rares prédécesseurs de Jordan à mettre en avant les "groupes de substitutions linéaires  définis par Galois" en relation avec les congruences et les équations cyclotomiques.\\
La référence  par Jordan à partir de 1861 à cette partie des travaux de Galois se présente donc comme une spécificité forte des travaux de ce dernier. Son approche générale  sur les groupes linéaires de $n$ variables n'allait d'ailleurs pas être reprise par la suite dans les présentations de la théorie des substitutions publiées jusqu'à la fin du siècle.\footnote{Voir notamment les présentations de Netto (1882), Klein (1884), Klein and Fricke (1890), Bolza (1891), Borel et Drach (1895), Weber (1895), Picard (1897), Pierpont (1900).} Comme nous l'évoquons dans le corps du texte, cette approche allait en réalité circuler dans un premier temps dans d'autres domaines comme les équations différentielles linéaires. 
\end{itemize}

\section{Annexe 4.\\
Sur le livre III Des irrationnelles \\ du \textit{Traité} de Jordan}
\label{Annexe4}

A partir de 1867, Jordan avait progressivement réorganisé ses travaux dans un nouveau cadre, centré non plus sur les relations transversales de la théorie de l'ordre mais sur un objet : le groupe de substitutions.

C'était dans ce nouveau cadre que se situait la présentation des principes généraux de Galois donnée au Livre III consacré aux irrationnelles. Cette approche se plaçait dans un héritage différent de celui de la théorie de l'ordre  qui consistait en une interprétation des "principes de Galois" en relation avec l'étude "irrationnelles" définies par des types d'équations algébriques générales non résolubles par radicaux (Hermite, Betti, Kronecker, Serret etc.). 

Contrairement à ce qui était devenu au XX\up{e} siècle un lieu commun de l'historiographie de l'algèbre, et à l'exception du domaine de l'enseignement de l'Algèbre supérieure, les travaux de Galois ont pendant longtemps été envisagés dans des cadres différents de ceux de la théorie des équations ou de la théorie des substitutions. Le problème de la "classification et la transformation" des "irrationnelles" impliquait notamment les fonctions elliptiques et abéliennes  et par conséquent l'analyse complexe. Il mêlait ainsi algèbre, arithmétique et analyse comme l'illustrent les résolutions par Hermite et Kronecker de l'équation générale du cinquième degré par l'équation modulaire d'ordre 5 des fonctions elliptiques. L'impossibilité d'exprimer par des fonctions algébriques les racines des équations de degré supérieur ou égal à 5 posait en effet la question de la détermination des fonctions analytiques les plus générales permettant l'expression de telles racines pour des types d'équations caractérisant des "ordres d'irrationalités".

Pour Jordan, l'invariance de la chaîne de quotients successifs de la décomposition des groupes par le théorème dit de Jordan-Hölder caractérisait l'"ordre d'irrationalité" d'une équation algébrique. Le Livre III s'appuyait alors à la fois sur l'approche d'Hermite mêlant analyse, algébrique et arithmétique\footnote{GOLDSTEIN (Catherine), SCHAPPACHER (Norbert), SCHWERMER (Joaquim) (eds.), \textit{The Shaping of Arithmetics after C. F. Gauss?s Disquisitiones Arithmeticae}, Berlin: Springer, 2007.} et sur l'approche géométrique de Clebsch sur la théorie des invariants et les fonctions abéliennes.

Le résultat clé du Livre III mettait en relation la réduction du degré de l'équation (de degré 80) de la trisection des périodes des fonctions abéliennes à quatre périodes et le groupe de l'équation des 27 droites sur une surface cubique. Cet énoncé jouait un rôle important pour légitimer l'approche de Jordan. Pour ce dernier, le résultat sur les 27 droites témoignait en effet de ce que les \textcolor{blue}{groupes de substitutions des équations spéciales de la géométrie}\footnote{\textcolor{blue}{Ces groupes n'étaient cependant pas introduits comme des groupes de Galois mais par l'invariance de formes algébriques. Le groupe symplectique $Sp_{2n}(p)$ était ainsi introduit par Jordan en 1869 comme laissant invariant une forme bilinéaire alternée non dégénérée. La notion d'"adjonction" de racine algébrique associée à Galois était quant à elle utilisée pour relier les groupes les uns aux autres. Par exemple, par "adjonctions" successives de racines, le groupe des 28 doubles tangentes d'une surface quartique se réduit aux groupes des 27 droites sur une cubique  et des 16 droites sur une quartique.}} permettent "l'étude des propriétés cachées de l'équation considérée".\cite[p.656]{16} 
\begin{quote}
Tous les géomètres connaissent le fait de l'abaissement des équations modulaires pour les transformations des degrés 5, 7 et 11, et les importantes conséquences qu'en a déduites M. Hermite. MM. Clebsch et Gordan ont signalé un abaissement analogue pour les équations des périodes dont dépend la bissection des fonctions abéliennes. Nous venons d'obtenir un résultat du même genre pour l'équation qui donne la trisection dans les fonctions à quatre périodes. [...] La réduite $Z$ présente cette particularité remarquable d'avoir le même groupe que l'équation $X$ qui détermine les vingt-sept droites situées sur une surface du troisième ordre. [...] Mais cette proposition exigeant quelques développements, nous en réservons la démonstration pour le \textit{Traité des équations algébriques} que nous nous occupons de publier.\cite{5}
\end{quote}
Jordan avait mené la réduction de l'équation des 27 droites sur le modèle de la réduction du degré des équations modulaires des fonctions elliptiques d'ordre 5. Comme l'avait énoncé Galois, et comme l'avait démontré Hermite à la suite de travaux de Betti, l'équation modulaire d'ordre 5 peut être "réduite" du degré 6 au degré 5. Une telle \textcolor{blue}{réduction}\footnote{\textcolor{blue}{En termes actuels, elle consiste à déterminer un sous groupe non distingué d'index $p$ de $PSl_2(p)$ (qui est simple pour $p>3$). Le groupe de Galois de l'équation modulaire réduite au degré 5 est $PSl_2(5)$, il est donc isomorphe au groupe alterné $A_5$ et peut être envisagé comme obtenu par une réduction du groupe symétrique $\Sigma(5)$ par l'"adjonction" du discriminant de l'équation générale du cinquième degré.}} est possible pour les équations modulaires d'ordres $p=5, 7, 11$. 

Une telle équation dont le groupe de Galois est plus petit que le groupe symétrique avait été dénommée "équation affectée" par Kronecker. La notion d'"affect" d'une équation avait été envisagée par ce dernier et par Hermite comme traduisant une propriété relative à l'ordre d'irrationalité de la quantité introduite par l'équation considérée. Cette notion visait initialement surtout les équations de degré 7. En effet, alors que l'équation générale de degré 5 peut être résolue par les fonctions elliptiques, des expressions analytiques de solutions d'équations de degré 7 ne pouvaient être obtenues que pour des cas d'équations affectées. L'un des principaux théorèmes du livre III de Jordan énonçait l'impossibilité de résoudre par les fonctions elliptiques les équations générales de degré strictement supérieur à 5 et la nécessité de recourir à des fonctions abéliennes de plusieurs périodes.

Par sa référence aux travaux de Galois, Jordan visait en partie à légitimer une approche qui, contrairement aux travaux d'Hermite, ne donnait pas plus de solution explicite au problème du nombre des valeurs des fonctions que d'expression analytique des racines d'une équation. Au contraire, Jordan revendiquait une approche "générale" de ces deux problèmes en considérant des substitutions de $n$ variables et en établissant des relations entre des classes de groupes (transitifs, primitifs, linéaires...). Ainsi, lors de son premier commentaire sur Galois dans le supplément de sa thèse, Jordan légitimait son premier théorème sur la décomposition des groupes primitifs en groupes linéaires sur $p^n$ lettres par son application aux équations algébriques: 
\begin{quote}
Il résulte de ce théorème que l'étude générale des équations et celle des systèmes de substitutions ne constituent au fond qu'un seul et même problème. Si les définitions et les théorèmes donnés dans ce Mémoire entrent bien dans le vif de la question, si les distinctions que j'y ai faites sont vraiment fondamentales, elles doivent représenter quelque propriété essentielle de l'équation correspondante au système considéré .\cite[p.187]{7}
\end{quote}
L'approche de Jordan avait été sévèrement critiquée par Kronecker comme une approche formelle procédant d'une fausse généralité. Pour ce dernier, Jordan avait donné abusivement à une méthode, l'usage de substitutions, le statut d'un objet d'étude.\footnote{Voir à ce sujet : BRECHENMACHER (Frédéric), La controverse de 1874 entre Camille Jordan et Leopold Kronecker, \textit{Revue d'Histoire des Mathématiques},  13 (2007), p. 187-257.}


\begin{thebibliography}{99}
\bibitem{10}
É.~Galois, Analyse d'un mémoire sur la résolution algébrique des équations, \textit{Bulletin des sciences mathématiques physiques et chimiques}, t. 13, n° 55 (1830), p. 171-172. Repr. in [Galois 1846, p.395-396].
\bibitem{11}
É.~Galois, Sur la théorie des nombres, \textit{Bulletin des sciences mathématiques physiques et chimiques}, t. 13, n° 55 (1830), p. 428-435. Repr. in [Galois 1846, p.398-407].
\bibitem{12}
É.~Galois, Fragment d'un second Mémoire. Des équations primitives qui sont solubles par radicaux, in [Galois 1846, p.434-444].
\bibitem{13}
É.~Galois,  Mémoire sur les conditions de résolubilité des équations par radicaux, in [Galois 1846, p. 417-433].
\bibitem{14}
É.~Galois, Lettre du 29 mai 1832 à Auguste Chevalier. \textit{Revue encyclopédique}, septembre 1832, p. 568-576. Repr. in [Galois 1846, p.  408-415]. 
\bibitem{15}
É.~Galois, \OE{}uvres  mathématiques, \textit{Journal de mathématiques pures et appliquées}, 11 (1846), 381-444.
\bibitem{1}
C.~Jordan, \textit{Sur le nombre des valeurs des fonctions, Thèses présentées à la Faculté des sciences de Paris par Camille Jordan, 1re thèse}, Paris : Mallet-Bachelier, 1860.
\bibitem{2}
C.~Jordan, Mémoire sur les groupes des équations solubles par radicaux, \textit{Comptes rendus hebdomadaires des séances de l'Académie des sciences}, t. 58 (1864), p. 963-966.
\bibitem{3}
C.~Jordan, Mémoire sur la résolution algébrique des équations, \textit{Journal de mathématiques pures et appliquées}, 12(2) (1867), pp. 109-157.
\bibitem{4}
C.~Jordan, Sur la résolution algébrique des équations primitives de degré $p^2$, \textit{Journal de mathématiques pures et appliquées}, 32 (2) (1868), p. 111-135.
\bibitem{16}
C.~Jordan, Sur les équations de la géométrie, \textit{Comptes rendus hebdomadaires des séances de l'Académie des sciences}, t.68 (1869), p. 656-659.
\bibitem{5}
C.~Jordan, Sur la trisection des fonctions abéliennes et sur les vingt-sept droites des surfaces du troisième ordre, \textit{Comptes rendus hebdomadaires des séances de l'Académie des sciences}, t.68 (1869), p. 865-869.
\bibitem{6}
C.~Jordan, \textit{Traité des substitutions et des équations algébriques}, Paris, 1870. 
\bibitem{7}
C.~Jordan, Sur la résolution des équations différentielles linéaires, \textit{Comptes rendus hebdomadaires des séances de l'Académie des sciences}, 73, 787-791.
\bibitem{8}
C.~Jordan, \textit{Notice sur les travaux de M. Camille Jordan à l'appui de sa candidature à l'Académie des sciences}, Paris: Gauthier-Villars, 1881.
\bibitem{9}
C.~Jordan, \textit{\OE{}uvres de Camille Jordan. Publiées sous la direction de M. Gaston Julia, par M. Jean Dieudonné}, Paris : Gauthier-Villars, 1961-1964.
\end{thebibliography}
\end{document}